\documentclass[10pt,a4paper]{amsart}
\usepackage{hyperref}
\usepackage[utf8x]{inputenc}
\usepackage{amsmath}
\usepackage{amssymb}
\usepackage{latexsym}
\usepackage{a4}
\usepackage[all]{xy}
\usepackage[T1]{fontenc}
\usepackage{mathrsfs}
\usepackage{eucal}
\title{On the double transfer and the $f$-invariant}
\author{Geoffrey Powell}
\address{Laboratoire Analyse, Géométrie et Applications, UMR 7539\\  Institut
Galilée, Université Paris 13, 93430 Villetaneuse, France}
\email{powell@math.univ-paris13.fr}
\urladdr{http://www.math.univ-paris13.fr/~powell/}
\thanks{Research partially supported by the project ANR BLAN08-2\_338236, HGRT}
\newtheorem{thm}{Theorem}[section]
\newtheorem{prop}[thm]{Proposition}
\newtheorem{cor}[thm]{Corollary}
\newtheorem{lem}[thm]{Lemma}

\theoremstyle{definition}
\newtheorem{defn}[thm]{Definition}
\newtheorem{exam}[thm]{Example}

\theoremstyle{remark}
\newtheorem{rem}[thm]{Remark}
\newtheorem{nota}[thm]{Notation}

\newcommand{\rhobar}{\overline{\rho}}
\newcommand{\qz}{{\rat/\zed}}
\newcommand{\qzsix}{{\rat/\zed[{\scriptstyle{\frac{1}{6}}}]}}
\newcommand{\zedsix}{{\zed[{\scriptstyle{\frac{1}{6}}}]}}
\newcommand{\zedpos}{\zed_{>0}}
\newcommand{\dash}{\hspace{-.2em}-\hspace{-.2em}}
\newcommand{\symm}{\mathfrak{S}}
\newcommand{\divBern}{\overline{B}}
\newcommand{\primel}{\mathfrak{p}}
\newcommand{\prim}{\mathbf{P}}
\newcommand{\divcong}{\mathfrak{D}}
\newcommand{\cobar}{\tilde{D}}
\newcommand{\BGammaC}{{_B\Gamma_C}}
\renewcommand{\epsilon}{\varepsilon}
\renewcommand{\phi}{\varphi}
\renewcommand{\hom}{\mathrm{Hom}}
\newcommand{\ext}{\mathrm{Ext}}
\newcommand{\modules}{\mathrm{Mod}}
\newcommand{\comod}{\mathrm{Comod}}
\newcommand{\zed}{\mathbb{Z}}
\newcommand{\cx}{\mathbb{C}}
\newcommand{\rat}{\mathbb{Q}}
\newcommand{\mf}{{MF}}
\newcommand{\mfmero}{\mf^{\mathrm{mer}}}
\newcommand{\padics}{\zed_p}
\newcommand{\plocal}{\zed_{(p)}}
\newcommand{\tmf}{{Ell}}
\newcommand{\tmfcoeff}{{Ell}}
\newcommand{\imj}{{Ad}}
\newcommand{\hq}{\mathbf{H}\mathbb{Q}}
\newcommand{\hz}{{\mathbf{H}\mathbb{Z}}}
\newcommand{\cp}{\mathbb{C} \mathbf{P}^\infty}

\newcommand{\gmu}{{\mathbb{G}_{\mathrm{m}, u}}}
\renewcommand{\smash}{\wedge}
\begin{document}

\begin{abstract}
The purpose of this paper is to investigate the algebraic  double $S^1$-transfer, in particular the classes in the two-line of the Adams-Novikov spectral sequence which are the image of comodule primitives of the $MU$-homology of $\cp\times \cp$ via the algebraic double transfer.

These classes are analysed by two related approaches; the first, $p$-locally for $p \geq 3$, by using the morphism induced in $MU$-homology by the chromatic factorization of the double transfer map together with the $f'$-invariant of Behrens (for $p \geq 5$). The second approach (after inverting $6$) uses the algebraic double transfer and the $f$-invariant of Laures.
\end{abstract}

\maketitle
\section{Introduction}

The double $S^1$-transfer is a stable morphism $\mathrm{tr}_2 : \cp _+ \smash \cp_+
\rightarrow S^{-2}$; determining its image in stable homotopy groups is a fundamental problem. This has an algebraic counterpart with respect to any complex oriented cohomology theory, in particular complex cobordism $MU$. Namely,  there is an algebraic double transfer $[e_\tau]^2$, which is a class in $\ext^2_{MU_* MU} (MU_* (\cp_+\smash \cp_+), MU_* [-4]) $, where $[-4]$ denotes the shift in internal degree and $\ext$ is calculated in the category of comodules over the Hopf algebroid $(MU_*, MU_ *MU)$;  this induces a morphism
\[
	\hom_{MU_* MU}^* (MU_*, MU_* (\cp_+ \smash  \cp_+))
	\rightarrow
	\ext^{2,*}_{MU_* MU} (MU_*, MU_* [-4])
\]
where the left hand side corresponds to the graded abelian group of $MU_*MU$-comodule primitives and the right hand side identifies with the $2$-line of the  Adams-Novikov $E^2$-term.

The corresponding algebraic framework for the single transfer is well understood, by the results of Miller \cite{miller}. For the double transfer, the situation is more complicated since the $MU_*MU$-comodule primitives of $MU_* (\cp_+ \smash \cp_+)$ are not fully understood and due to the additional complexity in passing from the Adams-Novikov $1$-line to the $2$-line. Baker approached this algebraic question in \cite{baker_transfer} by using Morava $K$-theory, working $p$-locally for a prime $p \geq 5$, in particular studying a family of primitives derived from the work of Knapp \cite{knapp_habilit}.

This paper approaches the algebraic question by two related methods, namely via the $f$-invariant of Laures \cite{laures} (which requires that $6$ is inverted) and via the $f'$-invariant of Behrens \cite{behrens}, which is defined when working $p$-locally for $p \geq 5$. Both are constructed by using elliptic homology, which is defined as a Landweber exact, complex oriented theory when $6$ is invertible. (In particular, for current purposes, it is not necessary to use topological modular forms.)

The $f'$-invariant  arises naturally when attempting to exploit the fact that the double transfer admits a chromatic factorization of the form
\[
	\cp_+ \smash \cp_+ \rightarrow S^{-4}/p^\infty, v_1 ^\infty
\]
for $p \geq 3$, which was first constructed by Hilditch (see \cite{BCGHRW}). The explicit determination of the induced morphism in $MU$-homology is non-trivial; it is determined implicity here by using  the Hattori-Stong theorem.

The first part of the paper explains this (see Theorem \ref{thm:MU-thetas} and Proposition \ref{prop:MUtheta_sigmap}) and relates it explicitly to calculations involving the algebraic double transfer, using standard chromatic technology. The result of Proposition \ref{prop:restrict_Xi_prim} is in principle sufficient to be able to calculate the algebraic double transfer on primitive elements; however, identifying the associated classes in $\ext^2$ is non-trivial (compare \cite{baker_transfer}).

In the second half of the paper, complex cobordism is replaced by elliptic homology and the $f$ and $f'$ invariants arising from comodule primitives of $\tmfcoeff_* (\cp_+ \smash \cp_+) $ are considered. The $f$-invariant on primitives is given by Theorem \ref{thm:f-pspt}, as a consequence of Proposition \ref{prop:restrict_Xi_prim}, whereas the $f'$-invariant on primitives is given (implicitly) by Theorem \ref{thm:fprime}. These results should shed light on the comodule primitives which are detected by the algebraic double transfer.

Behrens and Laures \cite{behrens_laures} have established the relationship between the $f$ and $f'$ invariants. For the invariants associated to comodule primitives via the double transfer, this relationship is transparent from chromatic technology, as indicated in Remark \ref{rem:relate_f,f'}.

\tableofcontents
\section{Chromatic factorization using $\mathrm{Im}(J)$}
\label{sect:hattori}

This section reviews the techniques for calculating morphisms to the spectrum $L_1 S /p^\infty$ of the
chromatic filtration, for $p$ a fixed prime,  and how to calculate the induced $MU_*MU$-comodule morphisms by  using the Hattori-Stong
theorem.

The terms ring spectrum and module spectrum refer to the weak, up to homotopy notions. If $E$ is
a ring spectrum and $M$ is an $E$-module the morphism of
$E$-module spectra induced by a morphism of spectra $f : X \rightarrow M$ is denoted $\tilde{f} : E \smash X \rightarrow M$.

\subsection{Non-connective $\mathrm{Im}(J)$-theory}

Let $\gamma \in \zed$ be  a topological generator of the $p$-adic
units $\padics^\times$. Non-connective image of $J$ theory, $\imj$, is
defined by the cofibre
sequence
\[
 \imj
\rightarrow
KU_{(p)}
\stackrel{\psi^\gamma -1} {\rightarrow}
KU_{(p)}
\rightarrow ,
\]
where $\psi^\gamma$ is the stable Adams operation, which is a morphism of ring spectra. The homotopy type of $\imj$ is independent of the choice
 of $\gamma$ (cf. \cite{knapp}).

The spectrum $\imj$ is a $KU$-module spectrum, in particular is
$KU$-local; moreover, there are  equivalences
\[
 \imj / p^\infty
 \simeq
\imj \smash S/p^\infty
\simeq
(L_1 S) \smash S/p^\infty
\simeq
L_1 (S/p^\infty)
\]
(cf. \cite[Lemma 8.7]{ravenel}), where  $L_
1$  is Bousfield localization with respect to $p$-local $K$-theory. Hence there
is a commutative diagram
\begin{eqnarray}
\label{eqn:Ad_cofibre_seq}
 \xymatrix{
&& L_1 S/p^\infty
\ar[d]^\alpha
\\
KU_{(p)}
\ar[r]
&
KU _\rat
\ar[r]^q
\ar[d]_{\psi^\gamma -1}
&
KU /p^\infty
\ar[d]^{\psi^\gamma - 1}
\\
&
KU_\rat
\ar[r]
&
KU /p^\infty
}
\end{eqnarray}
in which the three-term vertical and horizontal sequences are cofibre
sequences and $q$ is the reduction morphism. This provides a way of calculating maps to $L_1 S / p^\infty$,
as exploited in \cite[Section 5]{BCGHRW} and \cite{imaoka}, for example.

For the purposes of this paper, the following terminology is introduced.

\begin{defn}
	A $\rat$-representative of a morphism of spectra $g : Y \rightarrow L_1 S /p^\infty$ is a morphism $f : Y \rightarrow KU_\rat$ which makes the following diagram commute
	\begin{eqnarray}
\label{eqn:KU-commutative-ad}
 \xymatrix{
Y
\ar[r]^g
\ar[d]_f
&
L_1 S /p^\infty
\ar[d]^\alpha
\\
KU_\rat
\ar[r]_q
&
KU/p^\infty.
}
\end{eqnarray}
\end{defn}

\begin{lem}
	If $f$ is a $\rat$-representative of $g$, then $(\psi^\gamma -1) f$ lies in the image of $KU^0_{(p)}Y
\rightarrow KU^0_\rat Y$.
\end{lem}

\begin{proof}
	Follows from the commutativity of the square in diagram (\ref{eqn:Ad_cofibre_seq}).
\end{proof}

\begin{prop}
\label{prop:represent-ad}
Let $Y$ be a spectrum such that $KU ^*_{(p)} Y$ is a finitely-generated free $KU_{(p)*}$-module  and
$KU^{\mathrm{odd}}_{(p)}Y=0$. Then
\begin{enumerate}
	\item
	the morphism $[Y, L_1 S /p^\infty ] \rightarrow [Y, KU/p^\infty]$ is injective;
	\item
	any morphism $g : Y \rightarrow L_1 S /p^\infty$ admits a $\rat$-representative;
\item
a morphism $f: Y \rightarrow KU_\rat$ such that $(\psi^\gamma -1) f$ lies in the image of $KU^0_{(p)}Y \rightarrow
KU^0_\rat Y$ is the $\rat$-representative of a unique morphism $g : Y \rightarrow L_1 S/p^\infty$.
\end{enumerate}
\end{prop}

\begin{proof}
	Straightforward.
\end{proof}

\begin{exam}
	The hypotheses of Proposition \ref{prop:represent-ad} are satisfied for  $Y$ the Thom spectrum of  a finite rank virtual $\cx$-vector bundle over $\mathbb{C}\mathbf{P}^n$ and for smash products of spectra of this type.
\end{exam}

\subsection{Chromatic factorization}

Recall that a complex oriented ring spectrum $E$ is Landweber exact if the orientation $MU_* \rightarrow E_*$ is
Landweber exact for the Hopf algebroid $(MU_* , MU_*MU)$ (see Definition \ref{def:Landweber_exact}).

\begin{lem}
\label{lem:Landweber_exact}
	Let $E$ be a Landweber exact complex oriented ring spectrum and $Y$ be a spectrum.
	\begin{enumerate}
		\item There exist natural isomorphisms
		\[
			\hom_{MU_* MU} (MU_* Y, MU_* E) \cong \hom_{MU_*} (MU_* Y , E_*)
			\cong
			\hom_{E_*} (E_* Y, E_*).
		\]
\item
For a morphism of spectra $f : Y \rightarrow E$, the comodule morphism $MU_* f : MU_* Y \rightarrow MU_* E$ corresponds
via the above isomorphisms to the morphism of $E_*$-modules $\tilde{f}_* : E_* Y \rightarrow E_*$ induced by  $
\tilde{f}: E \smash Y \rightarrow E$.
\end{enumerate}
\end{lem}

\begin{proof}
The first isomorphism of part (1) follows from the identification of $MU_*E$ as the extended comodule $MU_*MU \otimes_{MU_*}  E_*$ and the second from the isomorphism of $E_*$-modules $E_* Y \cong E_* \otimes_{MU_*} MU_* Y$ which is a consequence of Landweber exactness. The final statement is straightforward.
\end{proof}

\begin{lem}
\label{lem:hattori-stong}
	Let $E$ be a Landweber exact complex oriented ring spectrum, then
the morphism $L_1 S /p^\infty \rightarrow KU/p^\infty$ induces a monomorphism of $E_* E$-comodules, $E_* /p^\infty [v_1^{-1}]  \hookrightarrow
E_* KU /p^\infty$.
\end{lem}

\begin{proof} By Landweber exactness,
	it suffices to prove this result for the universal case $E=MU$, where it is a consequence of the Hattori-Stong theorem (cf. \cite[Proposition 20.33]{switzer}), which states that the $KU$-Hurewicz morphism $MU_* \rightarrow MU_*KU$ is rationally faithful (in the terminology of \cite[Definition 1.1]{laures}), which is equivalent to the statement that $MU_* \otimes \rat/ \zed \hookrightarrow MU_* KU \otimes \rat/ \zed$ is a monomorphism. Hence, on the $p$-local component, this gives a monomorphism of $MU_*MU$-comodules $MU_* /p^\infty \hookrightarrow MU_* KU /p^\infty$.

	The morphism of $MU_*MU$-comodules $MU_* /p^\infty [v_1^{-1}] \rightarrow MU_* KU/p^\infty$ corresponds to the localization of the above morphism, inverting $v_1$, since the morphism $L_1 S \rightarrow KU_{(p)}$ factors the unit $S \rightarrow KU_{(p)}$. The result follows.
\end{proof}

\begin{prop}
\label{prop:f,g-Landweber}
Let $E$ be a Landweber exact complex oriented ring spectrum and $g: Y \rightarrow L_1 S /p^\infty$ be a morphism of spectra which admits a $\rat$-representative $f: Y \rightarrow KU_\rat$. Then
\begin{enumerate}
	\item
	the morphism
$E_* (g) : E_* Y \rightarrow E_* (L_1 S /p^\infty)$ is determined by $E_* (f)$ via the commutative diagram of morphisms of  $E_*E$-comodules:
\[
	\xymatrix{
E_* Y
\ar[r]^{E_* (g)}
\ar[d]_{E_*(f)}
&
E_*/p^\infty [v_1^{-1}]
\ar@{^(->}[d]
\\
E_* KU \otimes \rat
\ar@{->>}[r]
&
E_* KU/p^\infty .
}
\]
\item
The morphism $E_* (f)$  is determined by the morphism of $KU_*$-modules $\tilde{f}_* : KU_*  Y \rightarrow KU_*
\otimes \rat$ induced by $f$.
\end{enumerate}
\end{prop}

\begin{proof}
Again, by Landweber exactness, it  is sufficient to prove the result for the universal case, $E = MU$.

The commutative diagram (\ref{eqn:KU-commutative-ad}) induces a commutative diagram of $MU_* MU$-comodules:
\[
\xymatrix{
MU_* Y
\ar[r]^{MU_* (g)}
\ar@{^(->}[d]_{MU_*(f)}
&
MU_* /p^\infty [v_1 ^{-1}]
\ar@{^(->}[d]
\\
MU_*KU \otimes \rat
\ar@{->>}[r]
&
MU_* KU /p^\infty,
}
\]
in which $MU_* /p^\infty [v_1^{-1}] \rightarrow MU_* KU /p^\infty$ is
injective, by the Hattori-Stong theorem (Lemma \ref{lem:hattori-stong}),
 and $MU_* KU \otimes \rat \rightarrow
MU_* KU
/p^\infty$ is the canonical surjection. Thus, $MU_* (g)$ is determined by the
total composite of the diagram, hence by the morphism of $MU_*MU$-comodules
$MU_* Y \stackrel{MU_* (f)}{\rightarrow}MU_* KU \otimes
\rat$.

The final statement follows from Lemma \ref{lem:Landweber_exact}, which implies that the morphism $MU_*(f)$ is the composite
\[
	MU_* Y
	\stackrel{\psi}{\rightarrow}
	MU_* MU \otimes _{MU_*} MU_* Y
	\rightarrow
	MU_*MU \otimes _{MU_*} KU_* Y
	\stackrel{MU_* MU \otimes \tilde{f}_*}{\rightarrow}
	MU_*KU \otimes \rat,
\]
where $\psi$ is the comodule structure morphism and the second morphism is induced by $MU_* Y
\rightarrow KU_* Y$, given by the orientation of $KU$.
\end{proof}

In the case $E = KU$, this can be made more precise, by using the augmentation $KU_*KU \rightarrow KU_*$:

\begin{lem}
\label{lem:comodule-factorization}
Let $g: Y \rightarrow L_1 S /p^\infty$ be a morphism of spectra which admits a $\rat$-representative $f: Y \rightarrow KU_\rat$.
 Then there is an induced commutative diagram of
morphisms of $KU_*$-modules:
\[
 \xymatrix{
KU_* Y
\ar[r]^{KU_* (g)}
\ar[d]^{KU_* (f)}
\ar @{-->}@/_3pc/[dd]_{\tilde{f}_*}
&
KU_*/p^\infty
\ar[d]
\ar@/^3pc/[dd]^{KU_* /p^\infty}
\\
KU_*KU \otimes \rat
\ar[r]
\ar@{-->}[d]
&
KU_* KU/p^\infty
\ar@{-->}[d]
\\
KU_* \otimes \rat
\ar@{->>}[r]
&
KU_*/p^\infty
}
\]
in which the solid arrows are morphisms of $KU_*KU$-comodules and the lower
vertical morphisms are induced by the augmentation $KU_* KU \rightarrow KU_*$.

In particular, the comodule morphism $KU_* (g) : KU_* Y \rightarrow KU_*
/p^\infty$ factorizes as morphisms of $KU_*$-modules as
\[
 KU_ * Y
\stackrel{\tilde{f}_*} {\rightarrow}
KU_* \otimes \rat
\twoheadrightarrow
KU_* /p^\infty.
\]
\end{lem}

\section{Chromatic factorization of the double transfer}
\label{section:chromatic}

This section reviews the construction of the chromatic factorization of the double transfer (see Theorem
\ref{thm:chromatic-factor}), working $p$-locally at an odd prime $p$.  The morphism in $MU_*$-homology is
calculated implicitly  in $MU_*MU$-comodules (see Theorem \ref{thm:MU-thetas}), by applying the results of Section \ref{sect:hattori}.

\subsection{Generalized Bernoulli numbers}
\label{subsect:recollections}

To fix notation, recall that the Hopf algebroid $(MU_* , MU_*MU)$ is isomorphic to the Hopf algebroid $(L, LB)$ which represents the groupoid scheme
of formal group laws and strict isomorphisms, where $L$ is the Lazard ring and
$L B \cong L [b_i | i \geq 0, b_0 =1]$ as a left $L$-algebra (cf.
\cite{ravenel_green}). The $b_i's$ represent the universal strict isomorphism
$\underline{b} (x) = \sum_i b_i x^{i+1}$ between the formal group laws defined
respectively by the left and right units $\eta_L,\eta_R : L \rightrightarrows
LB$, which are determined by their logarithms $\log^L$, $\log^R$ defined over
$LB \otimes\rat$. The exponential series $\exp^L, \exp^R$ over $LB \otimes \rat$ are the respective
composition inverses of $\log^L, \log^R$.

\begin{lem}
\label{lem:identify_b}
The power series $\underline{b}$ satisfies the identity
$
\underline{b} = \exp^R \circ \log ^L.
$
\end{lem}

\begin{defn}
\cite[Definition 1.1]{miller}
Let $F$ be a formal group law defined over a ring $R$. The Bernoulli numbers
$B_n(F) \in R \otimes \rat$, for strictly positive integers $n \in \zedpos$,  are defined by
\[
 \frac{1}{\exp^F x} - \frac{1}{x}  =
\sum _{i\geq 0} \frac{B_{i+1}(F)}{(i+1)!} x^i,
\]
where $\exp^F(x) \in (R \otimes \rat) [[x]]$ is the exponential of $F$. The reduced Bernoulli number $\divBern_n (F)$  is defined as $\divBern_n (F) :=\frac{B_n (F)}{n} \in R \otimes \rat$.
\end{defn}

\begin{exam}
 For $n \in \zedpos$, write  $B_n^{KU} \in KU_* \otimes \rat$  (respectively
$\divBern_n ^{KU} \in KU_* \otimes \rat$)  for the  Bernoulli number (resp. reduced  Bernoulli number)  associated to the
orientation of  $KU$. This is a graded form of the usual  Bernoulli number $B_n$ (resp. reduced).
\end{exam}

\begin{rem}
	If the formal group law $F$ is graded with respect to the usual conventions (so that the coordinate has degree $-2$), then $B_n (F)$ is a homogeneous element of degree $2n$.
\end{rem}

\begin{rem}
	Miller  established the following fundamental divisibility property of the reduced Bernoulli numbers: if $R$ is a torsion free ring, then $d_n \divBern_n (F)  \in R$, where $d_n$ is the order of the reduced Bernoulli number $\divBern_n $ in $\qz$ (see \cite[Theorem 1.3]{miller}).
\end{rem}

\subsection{The single $S^1$-transfer}

\begin{defn}
For $n \in \zed$, let $\cp_n$ denote the Thom spectrum of the
(virtual) bundle $n \lambda$ over $\cp$, where $\lambda$ denotes the canonical
line bundle over $\cp$.
\end{defn}

For $E$  a complex oriented ring spectrum,  the Thom isomorphism implies that  $E_* (\cp_n)$ is a free $E_*$-module on
classes $\{\beta_i | i \geq n \}$.  (The systems of generators as $n$ varies are compatible, hence $n$ will be omitted from the notation.)
 There is a  Künneth  isomorphism $E_* (\cp_m \smash \cp _n) \cong E_* (\cp_m) \otimes _{E_*} E_* (\cp_n)$, and the associated module generators will be written $\beta_i \otimes \beta_j$.

\begin{nota}
	For $E$  a complex oriented ring spectrum and $m,n$  integers, let $\underline{\beta}_m (S)$ denote the
Laurent power series $\sum_{i \geq m} \beta_i S ^i$ over  $E_* (\cp_m)$ and let $\underline{\beta}_m (S) \otimes  \underline{\beta}_n (T)$ denote $\sum_{i \geq m, j  \geq n} \beta_i \otimes \beta_j S ^iT^j$, defined over $E_* (\cp_m \smash \cp_n)$.
\end{nota}

Such generating power series provide an efficient way of encoding calculations. For example:

\begin{lem}
\label{lem:comodule_cpn}
\cite[Proposition 3.3]{miller}
Let  $n$  be an integer, then the comodule structure $MU_* (\cp_n) \rightarrow MU_* MU \otimes_{MU_*} MU_*
(\cp_n)$ is determined by
\[
 \underline{\beta}_n (S)
\mapsto
\underline{\beta}_n (\underline{b} (S) \otimes 1).
\]
\end{lem}

\begin{rem}
In the expression $\underline{\beta}_n (\underline{b} (x) \otimes 1)$, the elements $\beta_i$ are  the
module generators, which are usually written on the right when considering  left $MU_*MU$-comodules.
Miller \cite{miller} works with right comodules, where this notational issue does not arise.
\end{rem}

The cofibre sequence of spectra (cf. \cite[section 2]{miller}):
\begin{eqnarray}
\label{eqn:cofibre-CPn}
S^{2n}
\rightarrow \cp_{n}
\rightarrow \cp_{n+1}
\rightarrow,
\end{eqnarray}
for $n \in \zed$, induces a short exact sequence of $MU_*MU$-comodules:
\[
0
\rightarrow
MU_* [2n]
\rightarrow
MU_* (\cp_n)
\rightarrow
MU_*(\cp_{n+1})
\rightarrow
0,
\]
where  $[a]$ denotes the shift in degree, so that $(V_* [a])_t = V_{t-a}$, for a $\zed$-graded object $V_*$.

The choice of generators  gives a standard splitting of this sequence in $MU_*$-modules. In particular:

\begin{nota}
\label{nota:sigma}
For $E$ a complex oriented ring spectrum,
\begin{enumerate}
	\item
let $\sigma : E_ * (\cp_0) \rightarrow E_* (\cp_{-1})$ be the section in $E_*$-modules
 defined by $\sigma (\beta_i) = \beta_i$ (for $i \geq 0$) and  $r : E_* (\cp_{-1}) \rightarrow E_*[-2]$
be  the corresponding retract, which  sends generators $\beta_i$, $i \geq 0$ to zero.
\item
let $\sigma' : E_ * (\cp_0 \smash \cp_0) \rightarrow E_* (\cp_{-1} \smash \cp_0)$ denote the section $\sigma \otimes E_* (\cp_0)$.
\end{enumerate}
\end{nota}

For $n=-1$, the connecting morphism of the cofibre sequence (\ref{eqn:cofibre-CPn}) defines  the $S^1$-transfer
$\tau : \cp_{0}\rightarrow S^{-1}$. The double $S^1$-transfer is
the smash product
$
	\tau \smash \tau :
\cp_{0} \smash \cp_0 \rightarrow S^{-2}.
$

The rational Thom class $U : \cp_{-1}\rightarrow S^{-2}_\rat$ induces a
morphism of cofibre sequences
\begin{eqnarray}
 \label{eqn:tau-tilde}
\xymatrix{
S^{-2}
\ar[r]
\ar@{=}[d]
&
\cp_{-1}
\ar[r]
\ar[d]_U
&
\cp_0
\ar[r]^\tau
\ar[d]^{\tilde{\tau}}
&
S^{-1}
\ar@{=}[d]
\\
S^{-2}
\ar[r]
&S^{-2}_\rat
\ar[r]
&
S^{-2}
_\qz
\ar[r]
&
S ^{-1},
}
\end{eqnarray}
where $\tilde{\tau}$, the chromatic factorization of the single transfer, is
determined uniquely by the commutativity of the right hand square.

The morphism   of $MU_*MU$-comodules
\[
	MU_* (\tilde{\tau})
:
MU_* (\cp_0)
\rightarrow
MU_* \otimes \qz [-2]
\]
is determined by the comodule morphism $MU_*(U) : MU_* (\cp_{-1})
\rightarrow
MU_* \otimes \rat [-2]$, via the commutative diagram
\begin{eqnarray}
\label{eqn:U-sigma}
	\xymatrix{
MU_* (\cp_{-1})
\ar[r]^{MU_*(U)}
&
MU_* \otimes \rat [-2]
\ar@{->>}[d]
\\
MU_* (\cp_0)
\ar@{.>}[u]^\sigma
\ar[r]_(.45){MU_* (\tilde{\tau}) }
&
MU_* \otimes \qz [-2],
}
\end{eqnarray}
in which the solid arrows denote  comodule morphisms.

By Lemma \ref{lem:Landweber_exact}, the comodule morphism $MU_*(U)$ is the composite
\[
\xymatrix{
	MU_* (\cp_{-1})
\ar[r]
&
MU_*MU \otimes _{MU_*} \hq_*(\cp_{-1})
\ar[d]^{MU_*MU \otimes \hq_*(U)}
\\
&
MU_*MU \otimes _{MU*} \rat [-2]
\ar[r]^(.6)\cong
&
MU_* \otimes \rat [-2],
}
\]
where $\hq$ is the rational Eilenberg-MacLane spectrum and the first morphism is the composite of the comodule structure morphism with  $MU_* (\cp_{-1} ) \rightarrow \hq_* (\cp_{-1})$ induced by the canonical orientation of $\hq$.

\begin{lem}(Cf. \cite[Theorem 3.9]{miller}.)
	\label{lem:MUThom}
	\begin{enumerate}
		\item
		The morphism $MU_* (U) :MU_* (\cp_{-1}) \rightarrow MU_* \otimes \rat [-2]$
is determined by
\[
	\underline{\beta}_{-1} (S)
\mapsto
\frac{1}{\log S}.
\]
\item
The morphism of $MU_*$-modules $MU_* (U) \circ \sigma : MU_* (\cp_0 ) \rightarrow MU_* \otimes \rat [-2]$
 is determined by
\[
	\underline{\beta}_{0} (S)
\mapsto
\frac{1}{\log S} - \frac{1}{S}.
\]
\item
The  morphism $MU_* (\tilde{\tau}) $ is the composite of $MU_*(U) \circ \sigma$ with the projection $MU_* \otimes \rat [-2]\twoheadrightarrow MU_* \otimes \qz[-2]$.
\end{enumerate}

\end{lem}

\begin{proof}
	 The first statement  follows from the comodule structure of $MU_*(\cp_{-1})$ together with the fact that,
under the morphism $MU_* MU \otimes \rat \cong MU_* \otimes MU_*
\otimes \rat \rightarrow MU_* \otimes \rat$ induced by the augmentation $MU_*
\otimes \rat \rightarrow \rat$ on the right hand factor, $\exp^R (S) \mapsto
S$, so that $\underline{b} (S) \mapsto \log (S)$.

The section $\sigma$ is determined by $\underline{\beta}_0 (S) \mapsto \underline{\beta}_{-1} (S) - \beta_{-1}\frac{1}{S}$, which gives the second statement, by composition. The final statement follows from the commutativity of diagram (\ref{eqn:U-sigma}).
\end{proof}

\subsection{The chromatic factorization of the double transfer}

Working $p$-locally ($p$ odd), the above chromatic factorization of the single transfer extends to a chromatic factorization of the double transfer, for which the original published reference is  \cite[Theorem 5.2]{BCGHRW}, where the result is attributed to Hilditch,  and a generalization is given by Imaoka in \cite{imaoka}. (Imaoka \cite{imaoka_p2} has also considered the chromatic factorization of the double transfer at the prime $p=2$.)

\begin{thm}
\cite{BCGHRW}
\label{thm:chromatic-factor}
Let $p \geq 3$ be a prime.
There exists a morphism $\Theta : \cp_{-1} \smash \cp_0 \rightarrow L_1 S^{-4}/p^\infty$ which fits into a commutative square:
\begin{eqnarray}
\label{eqn:Theta_extends_tautilde}
\xymatrix{
\Sigma^{-2} \cp_0
\ar[r]
\ar[d]_{\Sigma^{-2}\tilde{\tau}}
&
\cp_{-1} \smash \cp_0
\ar[d]^\Theta
\\
S^{-4} /p^\infty
\ar[r]
&
L_1 S ^{-4}/p^\infty,
}
\end{eqnarray}
where the top morphism is induced by the inclusion of the bottom cell $S^{-2} \rightarrow \cp_{-1}$.

Moreover, for any extension to a  morphism of cofibre sequences:
\[
\xymatrix{
\Sigma^{-2} \cp_0
\ar[r]
\ar[d]_{\Sigma^{-2}\tilde{\tau}}
&
\cp_{-1} \smash \cp_0
\ar[d]^\Theta
\ar[r]
&
\cp_0 \smash \cp_0
\ar[d]^{\overline{\Theta}}
\ar[r]^(.6){\tau \smash \cp_0}
&
\Sigma ^{-1}\cp_0
\ar[d]^{\Sigma^{-1} \tilde{\tau}}
\\
S^{-4} /p^\infty
\ar[r]
&
L_1 S ^{-4}/p^\infty
\ar[r]
&
S^{-4} /p^\infty , v_1^\infty
\ar[r]
&
 S^{-3}/p^\infty,
}
\]
where the top row is the cofibre sequence $(S^{-2}\rightarrow
\cp_{-1} \rightarrow \cp_0)\smash \cp_0 $ and the bottom row is the cofibre
sequence associated to $L_1$-localization,  $\overline{\Theta} :
\cp_0 \smash \cp_0
\rightarrow S^{-4} /p^\infty , v_1^\infty$ provides a factorization of
the double transfer morphism across the chromatic morphism $S^{-4}
/p^\infty , v_1^\infty \rightarrow S^{-2}$.
\end{thm}

The morphism $\Theta$ is constructed using Proposition \ref{prop:represent-ad}, by defining an
explicit cohomology class cohomology class $\theta \in [\cp_{-1} \smash \cp_0 , \Sigma^{-4} KU _\rat ]$ such that there
is a commutative diagram
\begin{eqnarray}
\label{eqn:theta_Q_representative}
\xymatrix{
\cp_{-1}\smash \cp_0
\ar[r]^\Theta
\ar[d]_\theta
&
L_1 S ^{-4}/p^\infty
\ar[d]
\\
\Sigma^{-4} KU _\rat
\ar[r]
&
\Sigma^{-4} KU /p^\infty.
}
\end{eqnarray}

\begin{rem}
	The construction of $\theta$ is by an eigenspace argument for the action of the Adams operation $\psi^\gamma$,
where $\gamma \in \zed$ is a topological generator of $\padics$ (compare  Proposition \ref{prop:represent-ad}); this is made explicit in the proof of \cite[Proposition 2.4]{imaoka}, which generalizes this result.
\end{rem}

For later use, the following notation is introduced.

\begin{nota}
\label{nota:tilde_theta_prime}
	Let $\tilde{\theta} ' (S,T)  $ denote the power series in $KU_* \otimes \rat [[S,T]]$
	\[
 \sum _{i, j >0}\frac{B_i^{KU}}{i!}\frac{B_j^{KU}}{j!}
\left(
\frac{\gamma^i -1}{\gamma^{i+j}-1}
\right)
(\log^{KU} S)^{i-1}(\log ^{KU} T) ^{j-1},
\]
where $\log^{KU}$ is the logarithm of the multiplicative formal group law of $KU_* \otimes \rat$.
\end{nota}

The morphism $\theta$ is not uniquely defined; \cite[Theorem 5.2]{BCGHRW} gives an explicit choice for $\theta$, working
with the $p$-local Adams summand $E(1)$, which can be replaced by $p$-local $K$-theory. The following choice is used here:

\begin{defn}
\label{def:theta}
Let $\theta : \cp_{-1} \smash \cp_0 \rightarrow \Sigma ^{-4} KU_\rat$ be the class
 which is determined by
$\tilde{\theta}_* : KU_*(\cp_{-1}\smash \cp_0) \rightarrow KU_* \otimes \rat [-4]$:
\[
	\underline{\beta}_{-1}(S) \otimes \underline{ \beta}_0(T)
	\mapsto
	\frac{1}{S} \left(\frac{1}{\log^{KU}T} - \frac{1}{T}\right)
+
\tilde{\theta}' (S,T).
\]
\end{defn}

\subsection{Calculating in complex cobordism}

Let $p$ be a fixed odd prime and $\Theta, \overline{\Theta}$ be as in  Theorem
\ref{thm:chromatic-factor}, where  the $\rat$-representative $\theta$ of $\Theta$ is the morphism of Definition
\ref{def:theta}.

\begin{thm}
\label{thm:MU-thetas}
Let $p$ be an odd prime.
There is  a commutative diagram of morphisms of
$MU_*MU$-comodules:
\[
\xymatrix{
 &
MU_* (\cp_{-1}\smash \cp_0)
\ar@{->>}[r]
\ar@/_1pc/[ld]_(.6){MU_*(\theta)}
\ar[d]^{MU_*(\Theta)}
&
MU_* (\cp_{0}\smash \cp_0)
\ar[d]^{MU_*(\overline{\Theta})}
\\
MU_* KU \otimes \rat[-4]
\ar@{->>}@/_1pc/[rd]
&
MU_*/p^\infty[v_1^{-1}] [-4]
\ar@{^(->}[d]
\ar@{->>}[r]^{\mathrm{pr}}
&
MU_* /p^\infty,v_1^\infty [-4]
\ar@{^(->}[d]
\\
&
MU_* KU /p^\infty[-4]
\ar@{->>}[r]
&
(MU_* KU /p^\infty)/(MU_*/p^\infty)[-4].
}
\]
\begin{enumerate}
	\item
	The underlying $MU_*$-module morphism of  $MU_* (\overline{\Theta})$ is the composite $\mathrm{pr}\circ MU_* (\Theta)
\circ \sigma'$.
\item
\label{thm_item:Theta}
The morphism $MU_*({\Theta})$ is determined by the comodule
morphism $MU_* (\theta)$ and hence by the morphism $\tilde{\theta}_* : KU_* (\cp_{-1}\smash \cp_0
)\rightarrow KU_* \otimes \rat$.
\item
The morphism $MU_*(\overline{\Theta})$ is determined by the comodule
morphism $MU_* (\theta)$.
\end{enumerate}
\end{thm}

\begin{proof}
	The commutativity of the left hand part of the diagram follows from Proposition \ref{prop:f,g-Landweber} and the  upper right hand square is induced by the morphism of cofibre sequences defining $\overline{\Theta}$.  The lower right hand square is induced by the monomorphism given by the Hattori-Stong theorem (see Lemma \ref{lem:hattori-stong})
\[
	MU_* /p^\infty [v_1^{-1}] \hookrightarrow MU_ * KU /p^\infty,
\]
since the kernel to the projection $MU_* /p^\infty [v_1^{-1}] \twoheadrightarrow MU_*/p^\infty , v_1^\infty $ is
$MU_* /p^\infty$.

The morphism $MU_* (\cp_{-1} \smash \cp_0) \twoheadrightarrow MU_* (\cp_0 \smash \cp_0)$  admits the section $\sigma'$ in $MU_*$-modules, hence the identification of the underlying module morphism of $MU_* (\overline{\Theta})$ follows from the upper right hand square.

For part \ref{thm_item:Theta}, the injectivity of
$MU_* /p^\infty, v_1 ^\infty [-4] \hookrightarrow MU_* KU /p^\infty[-4]$
and the commutativity of the left hand square implies that $
MU_* (\theta)$ determines $MU_* ({\Theta})$ as a morphism of $MU_*MU$-comodules. The morphism $MU_* (\theta)$ is determined by $\tilde{\theta}_* : KU_* (\cp_{-1}\smash \cp_0
)\rightarrow KU_* \otimes \rat$, by the second part of Proposition \ref{prop:f,g-Landweber}.

The final statement follows from the commutativity of the lower right hand square.
\end{proof}

\subsection{Calculating $MU_* (\theta)$}

The torsion-free ring $MU_*KU$ has two formal group law structures, corresponding to the left and right units $MU_* \rightarrow
MU_*KU$ and $KU_* \rightarrow MU_* KU$ respectively; write $\log^{MU}$ and $\log^{KU}$ for the respective  logarithms
defined  over $MU_* KU \otimes \rat$ and set  $\underline{b}' = \exp^{KU} \circ \log^{MU}$, which identifies
with the image of $\underline{b}$ under $MU_*MU \rightarrow MU_*KU$.

\begin{prop}
\label{prop:MUtheta_sigmap}
\
\begin{enumerate}
\item
The morphism $MU_*(\theta)\in \hom_{MU_*MU} (MU_* (\cp_{-1} \smash \cp_0 ), MU_* KU \otimes \rat[-4]) $ is
determined by
$
 \underline{\beta}_{-1} (S) \otimes \underline{\beta}_0 (T)
\mapsto $
\[
 \frac{1}{\underline{b}'(S)}
\Big(\frac{1}{\log^{MU} T} - \frac{1}{\underline{b}' (T) }\Big)
+
\sum _{i, j >0}\frac{B_i^{KU}}{i!}\frac{B_j^{KU}}{j!}
\left(
\frac{\gamma^i -1}{\gamma^{i+j}-1}
\right)
(\log^{MU} S)^{i-1}(\log ^{MU} T) ^{j-1}.
\]
 \item
The morphism $MU_*(\theta) \circ \sigma' \in \hom_{MU_*} (MU_* (\cp_0 \smash \cp_0 ), MU_* KU \otimes \rat[-4]) $ is
determined by
$
 \underline{\beta}_0 (S) \otimes \underline{\beta}_0 (T)
\mapsto $
\[
\Big( \frac{1}{\underline{b}'(S)}- \frac{1}{S}\Big)
\Big(\frac{1}{\log^{MU} T} - \frac{1}{\underline{b}' (T) }\Big)
+
\sum _{i, j >0}\frac{B_i^{KU}}{i!}\frac{B_j^{KU}}{j!}
\left(
\frac{\gamma^i -1}{\gamma^{i+j}-1}
\right)
(\log^{MU} S)^{i-1}(\log ^{MU} T) ^{j-1}.
\]
\end{enumerate}
\end{prop}

\begin{proof}
 By Proposition \ref{prop:f,g-Landweber}, the morphism $MU_* (\theta)$
is determined by the commutative diagram

\[
\xymatrix{
	MU_* (\cp_{-1} \smash \cp_0)
\ar[rr]^(.4){\psi _{MU_* (\cp_{-1} \smash \cp_0)}}
\ar[d]_{MU_* (\theta)}
&\ &
MU_*MU
\otimes_{MU_*} MU_*(\cp_{-1} \smash \cp_0)
\ar[d]
\\
MU_*MU\otimes_{MU_*} KU_* \otimes \rat [-4]
&&
MU_*MU \otimes _{MU_*} KU_* (\cp_{-1} \smash \cp_0)
\ar[ll]^{MU_*MU \otimes \tilde{\theta}_*}.
}
\]
Lemma \ref{lem:comodule_cpn} implies that the comodule structure of $MU_* (\cp_{-1} \smash \cp_0)$ is determined by
\[
	\underline{\beta}_{-1} (S) \otimes \underline{\beta}_0 (T)
	\mapsto
	\underline{\beta}_{-1} (\underline{b}(S)\otimes 1) \otimes \underline{\beta}_0 (\underline{b}(T) \otimes 1).
\]
Composing with the morphism induced by $\theta$ and using the identity $\underline{b}' = \exp^{KU} \circ \log^{MU}$ shows that $MU_* (\theta)$ is given by
$\underline{\beta}_{-1} (S) \otimes \underline{\beta}_0 (T)
	\mapsto$
\[
\frac{1}{\underline{b}'(S)}
\Big(\frac{1}{\log^{MU} T} - \frac{1}{\underline{b}' (T) }\Big)
+
\sum _{i, j >0}\frac{B_i^{KU}}{i!}\frac{B_j^{KU}}{j!}
\left(
\frac{\gamma^i -1}{\gamma^{i+j}-1}
\right)
(\log^{MU} S)^{i-1}(\log ^{MU} T) ^{j-1}.
\]
The second statement is proved by composing with the  morphism $\sigma'$, which  is represented by
\[
	\underline{\beta}_0 (S) \otimes \underline{\beta}_0 (T) \mapsto \underline{\beta}_{-1} (S) \otimes \underline{\beta}_0 (T) - \beta_{-1} \frac{1}{S}
	\otimes \underline{\beta}_0 (T).
\]
The morphism $MU_* (\theta)$ restricts to give:
\[
	\beta_{-1} \otimes \underline{\beta}_0 (T) \mapsto \frac{1} {log^{MU} T} - \frac{1}{\underline{b}'(T)}.
\]
The result follows.
\end{proof}

\section{The algebraic transfer}
\label{sect:cohomtransfer}

The algebraic version of the double transfer is introduced in this section and is related to chromatic theory.

\subsection{The algebraic transfer}
\label{subsect:alg_transfer}

The cofibre sequence defining the transfer $\tau : \cp_0 \rightarrow
S^{-1}$ induces a short exact sequence of $MU_*MU$-comodules and hence an
algebraic transfer class:
\[
[e_\tau] \in \ext^1_{MU_* MU}(MU_* (\cp_0),MU_* [-2]).
\]

\begin{defn}
\label{def:algebraic_double_transfer}
The algebraic double transfer is the class:
\[
 [e_\tau]^2\in \ext^2_{MU_* MU}(MU_* (\cp_0\smash \cp_0 ),MU_* [-4])
\]
given by Yoneda product.
\end{defn}

Proposition \ref{prop:cocycle-canonical-splitting} applied with respect
to the section $\sigma$ gives the standard choice $e_\tau$ of representing cocycle:

\begin{lem}
\label{lem:cocycle-etau}
The cocycle $
{e_\tau} \in \hom_{MU_*}(MU_* (\cp_0), MU_* MU[-2])
$
is determined by
\[
	\underline{\beta}_{0} (S)
\mapsto
\frac{1}{S}
-
\frac{1}{\underline{b} (S)} .
\]
\end{lem}

Diagram (\ref{eqn:tau-tilde}) induces a morphism of short exact sequences of $MU_*MU$-comodules:
\[
	\xymatrix{
0
\ar[r]
&
MU_* [-2]
\ar@{=}[d]
\ar[r]
&
MU_* (\cp_{-1} )
\ar[d]^{MU_* (U)}
\ar[r]
&
MU_* (\cp_0)
\ar[d]^{MU_*(\tilde{\tau})}
\ar[r]
&
0
\\
0
\ar[r]
&
MU_* [-2]
\ar[r]
&
MU_* \otimes \rat [-2]
\ar[r]
&
MU_* \otimes \qz [-2]
\ar[r]
&
0.
}
\]
The morphism $\tilde{\tau}$ provides a chromatic factorization of the single transfer $\tau$; this corresponds to the
following result, which can be proved using Proposition \ref{prop:cocycle-canonical-splitting}.

\begin{prop}
	There is an equality in $\ext^1_{MU_* MU} (MU_* (\cp_0), MU_* [-2])$:
\[
	[e_\tau] = \partial_1 MU_* (\tilde{\tau}),
\]
where $\partial _1$ is the  chromatic connecting morphism
associated to
\[
\xymatrix{
0
\ar[r]
&
MU_* [-2]
\ar[r]
&
MU_* \otimes \rat [-2]
\ar[r]
&
MU_* \otimes \qz [-2]
\ar[r]
&
0.
}
\]
\end{prop}

\subsection{The class $[\kappa]$}

Rather than working directly with the double algebraic transfer $[e_\tau]^2$, it is convenient to
work with a  class $[\kappa]$ in $\ext^1$ (see Definition \ref{def:kappa}), which is related to
the double transfer  via the chromatic connecting morphism $\partial_1$ (see Proposition \ref{prop:kappa-transfer}).

Forming the tensor product $[e_\tau] \otimes MU_* (\cp_0)$ gives a
class
\[
[E_\tau ]\in \ext^1_{MU_* MU}(MU_* (\cp_0 \smash \cp_0),MU_*(\cp_0) [-2]).
\]

\begin{lem}
\label{lem:Etau}
The class $[E_\tau]$ is represented by the cocycle
\[
{E_\tau} \in \hom_{MU_*}(MU_* (\cp_0\smash \cp_0), MU_* MU\otimes _{MU_*} MU_* (\cp_0)[-2])
\]
defined with respect to the section $\sigma'$,  which is given by
\[
 \underline{\beta}_ 0 (S) \otimes \underline{\beta}_0 (T)
\mapsto
\Big(
\frac{1}{S}
- \frac{1}{\underline{b}(S)}
\Big)
\underline{\beta}_0 (\underline{b}(T)).
\]
\end{lem}

\begin{proof}
Apply Proposition \ref{prop:cocycle-canonical-splitting} with respect to the section $\sigma'$.
\end{proof}

\begin{defn}
\label{def:kappa}
	Let $[\kappa]$ denote the class
\[
MU_* (\tilde{\tau}) [E_\tau] \in \ext^1
_{MU_* MU} (MU_* (\cp_0 \smash \cp_0), MU_* \otimes \qz[-4]).
\]
\end{defn}

The following result justifies using $\kappa$ in place of the algebraic double transfer.

\begin{prop}
 \label{prop:kappa-transfer}
There is an identity
\[
 [e_\tau]^2 = \partial_1 [\kappa],
\]
in $\ext^2 _{MU_*MU}(MU_* (\cp_0 \smash \cp_0), MU_* [-4])$.
\end{prop}

\begin{proof}
 Straightforward.
\end{proof}

The class $[\kappa]$ is represented by the standard choice  $\kappa$ of cocycle, constructed with respect
to the section $\sigma'$:
\[
{\kappa} \in \hom_{MU_*} ( MU_* (\cp_0\smash \cp_0),MU_*MU \otimes \qz [-4]).
\]

\begin{nota}
\label{nota:Xi}
	Write $K \in \hom_{MU_*} (MU_* (\cp_0\smash \cp_0), MU_*MU \otimes
\rat [-4])$ for the composite  morphism  of $MU_*$-modules:
\[
	\xymatrix{
MU_* (\cp_0 \smash \cp_0)
\ar[r]^(.4){E_\tau}
\ar[rd]_K
&
MU_* MU\otimes _{MU_*}MU_*(\cp_{0})[-2]
\ar[d]^{MU_*MU \otimes MU_*(U) \circ \sigma }
\\
&
MU_*MU \otimes \rat [-4].
}
\]
\end{nota}

\begin{prop}
 \label{prop:calculate_Xi}
\
\begin{enumerate}
	\item
The class $[\kappa]$ is represented by the cocyle $\kappa$
given by reduction  of the morphism $K$ via $MU_*MU \otimes \rat  [-4]
\rightarrow
MU_*MU \otimes \qz[-4]$.
\item
The morphism $K \in \hom_{MU_*} (MU_* (\cp_0\smash \cp_0), MU_*MU \otimes
\rat [-4])$ is determined by
\[
 \underline{\beta}_0 (S) \otimes \underline{\beta}_0 (T)
\mapsto
\Big(
\frac{1}{S}- \frac{1}{\underline{b}(S)}
\Big)
\Big(
\frac{1}{\log^L T }
- \frac{1}{\underline{b}(T)}
\Big).
\]
\end{enumerate}
\end{prop}

\begin{proof}
The first statement follows from the construction of $K$.

The second follows from Lemma \ref{lem:Etau}, by composition with the morphism
$MU_*MU \otimes (MU_*(U) \circ\sigma)$ (see Lemma \ref{lem:MUThom}), using the identity $\log^L = \log^R \circ \underline{b}$ of Lemma \ref{lem:identify_b}.
\end{proof}

\subsection{Relation with the chromatic factorization of the double transfer}
\label{subsect:relate-double-chromatic}

Let $p$ be an odd prime. Here $\kappa$ will be written to denote the associated $p$-local cocycle
\[
	\kappa
\in \hom_{MU_*} (MU_* (\cp_0\smash \cp_0) , MU_*MU/p^{\infty}[-4]).
\]

The construction of $\overline{\Theta}$ as a chromatic factorization
of the double transfer implies that the  class $[\kappa]$ is related to the morphism $MU_*(\overline{\Theta})$. Write
$\partial_2$ for the chromatic connecting morphism associated to the short exact
sequence of comodules
\[
	0
\rightarrow MU_* /p^\infty
 \rightarrow
MU_* /p^\infty [v_1 ^{-1}]
\rightarrow
MU_* /p^\infty ,v_1 ^\infty
\rightarrow 0.
\]

\begin{prop}
\label{prop:second-factorization}
	There is an identity
$
	[\kappa ] = \partial _2 MU_* (\overline{\Theta})
$
in $\ext^1_{MU_*MU} (MU_* (\cp_0\smash \cp_0 ), MU_*/p^\infty [-4])$. In particular, $\partial_2 MU_* (\overline{\Theta}) $ is independent of the choice of $\overline{\Theta}$.
\end{prop}

\begin{proof}
The morphism $MU_* (\Theta)$ gives rise to a morphism between short exact sequences of $MU_*MU$-comodules:
\[
\xymatrix{
0
\ar[r]
&
MU_* (\cp_0)[-2]
\ar[d]_{MU_* (\tilde{\tau}) [-2]}
\ar[r]
&
MU_* (\cp_{-1} \smash \cp_0)
\ar[r]
\ar[d]_{MU_* (\Theta)}
&
MU_* (\cp_0 \smash \cp_0)
\ar[d]^{MU_* (\overline{\Theta}) }
\ar[r]
&
0
\\
0
\ar[r]
&
MU_*/p^\infty[-4]
\ar[r]
&
MU_* /p^\infty [v_1]^{-1}[-4]
\ar[r]
&
MU_* /p^\infty , v_1 ^\infty [-4]
\ar[r]
&
0,
}
\]
where the top row represents $[E_\tau]$. By definition, $[\kappa] = MU_* (\tilde{\tau}) [E_\tau]$ and $\partial _2 MU_* (\overline{\Theta})$ is represented by the pullback of the lower short exact sequence along $MU_* (\overline{\Theta})$.

Forming the pushout of the top sequence using $MU_* (\tilde{\tau}) $ and the pullback of the lower
sequence via $MU_* (\overline{\Theta})$ gives Yoneda-equivalent short exact sequences,  which therefore define the same
class in $\ext^1$, as required.
\end{proof}

\begin{rem}
	It is instructive to check this result directly at the level of cocycles by using the description of the connecting morphism given in Lemma \ref{lem:explicit-connecting}.
\end{rem}

\subsection{Relating $K$ and $MU_*(\theta) \circ \sigma'$}
\label{subsect:relate_K_MuTheta}

The morphism $MU_*MU \rightarrow MU_* KU$ associated to the
orientation of $KU$ induces the morphism (using the notation introduced in Section \ref{sect:cocycle-basechange}):
\[
 _{MU_*} K _{KU_*} \in \hom_{MU_*} (MU_* (\cp_0\smash \cp_0), MU_*KU \otimes
\rat [-4]).
\]
This can be identified with the composite
\[
 \xymatrix{
MU_* (\cp_0 \smash \cp_0)
\ar[r]^{- \sigma'}
&
MU_* (\cp_{-1} \smash \cp_0)
\ar[r]^(.35){\psi}
&
MU_*MU \otimes_{MU_*}
MU_* (\cp_{-1} \smash \cp_0)
\ar[d]
\\
&
MU_* KU \otimes \rat[-4]
&
MU_*MU \otimes_{MU_*}
KU_* (\cp_{-1} \smash \cp_0)
\ar[l],
}
\]
where $\psi$ is the comodule structure map, the vertical arrow is induced by the orientation of $KU$ and the final
morphism of $MU_*MU$-modules is defined by the composite
\begin{eqnarray}
\label{eqn:KU_thom_sigma}
 KU_* (\cp_{-1} \smash \cp_0)
\twoheadrightarrow
KU_* (\cp_0) [-2]
\stackrel{KU_*(U) \circ \sigma}{\longrightarrow}
KU_* \otimes \rat [-4]
\end{eqnarray}
of the projection induced by $KU_* (\cp_{-1}) \twoheadrightarrow KU_* [-2]$ with the morphism induced by the rational
Thom class of $\cp_{-1}$.

\begin{rem}
The sign arises due to the conventions used in defining the cobar complex, as in Proposition \ref{prop:cocycle-canonical-splitting}.
\end{rem}

The second  morphism in (\ref{eqn:KU_thom_sigma}) is related to the morphism $\tilde{\theta}_* : KU_* (\cp_{-1} \smash \cp_0) \rightarrow KU_* \otimes
\rat [-4]$ via the following commutative diagram derived from Theorem \ref{thm:MU-thetas}:
\[
 \xymatrix{
\Sigma^{-2} \cp_{-1}
\ar[r]
\ar[d]_{\Sigma^{-2}U}
&
\Sigma^{-2} \cp_0
\ar[d]_{\Sigma^{-2}\tilde{\tau}}
\ar[r]
&
\cp_{-1}
\smash
\cp_0
\ar[d]^{\Theta}
\ar[ldd]|\hole_(.3)\theta
\\
S^{-4}_\rat
\ar[r]
\ar[rd]
&
S^{-4}/p^\infty
\ar[r]
&
L_1 S^{-4} /p^\infty
\ar[d]
\\
&
\Sigma^{-4} KU_\rat
\ar[r]
&
\Sigma^{-4} KU/p^\infty.
}
\]

This implies the following result (which corresponds to a fundamental property of $\theta$ used in the construction of $\Theta$).

\begin{lem}
\label{lem:restrict_Theta}
The restriction of $\tilde{\theta}_* : KU_* (\cp_{-1} \smash \cp_0) \rightarrow
KU_* \otimes \rat[-4]$ along the morphism $KU_* (\cp_0) [-2] \hookrightarrow KU_* (\cp_{-1} \smash \cp_0)$, induced by
the inclusion of the bottom cell $S^{-2} \hookrightarrow \cp_{-1}$ is
 the morphism
 \[
 KU_*(U) \circ \sigma[-2] : 	KU_* (\cp_0) [-2]
\rightarrow KU_* \otimes \rat [-4].
 \]
\end{lem}

\begin{rem}
This result corresponds to the fact that the morphism $\tilde{\theta}_*$ is determined by the  power series
	\[
		\frac{1}{S} \left(\frac{1}{\log^{KU}T} - \frac{1}{T}\right)
+
\tilde{\theta} ' (S,T)
	\]
where $\tilde{\theta}'$ is the formal power series introduced in Notation \ref{nota:tilde_theta_prime}.
\end{rem}

Proposition \ref{prop:calculate_Xi} gives the following, using the notation of Proposition \ref{prop:MUtheta_sigmap}:

\begin{prop}
 \label{prop:calculate_K_KU}
The morphism $$_{MU_*}K_{KU_*} \in \hom_{MU_*} (MU_* (\cp_0\smash \cp_0), MU_*KU \otimes
\rat [-4])$$ is determined by
\[
 \underline{\beta}_0 (S) \otimes \underline{\beta}_0 (T)
\mapsto
\Big(
\frac{1}{S}- \frac{1}{\underline{b}'(S)}
\Big)
\Big(
\frac{1}{\log^{MU} T }
- \frac{1}{\underline{b}'(T)}
\Big).
\]
\end{prop}

\begin{proof}
 By construction $\underline{b}'$ corresponds to the image of $\underline{b}$ under the morphism induced by $MU_*MU
\rightarrow MU_*KU$, and $\log^L$ maps to $\log^{MU}$. The result follows  from Proposition
\ref{prop:calculate_Xi}.
\end{proof}

The above description of $ _{MU_*} K _{KU_*}$ can be compared with that of $MU_* (\theta ) \circ \sigma ' : MU_* (\cp_0
\smash \cp_0) \rightarrow MU_*KU \otimes \rat [-4]$. Write $\tilde{\theta}'_*$ for the  morphism of $MU_*$-modules $MU_* (\cp_0 \smash \cp_0) \rightarrow KU_* \otimes \rat[-4]$ determined by  $\tilde{\theta}'$.

\begin{cor}
\label{cor:relate_MUtheta}
There is an identification of morphisms in $\hom_{MU_*} (MU_* (\cp_0\smash \cp_0), MU_*KU \otimes
\rat [-4])$:
\[
	MU_* (\theta) \circ \sigma'
	=
	-\  _{MU_*}K_{KU_*} + (MU_*MU \otimes \tilde{\theta}'_*) \circ \psi_{MU_* (\cp_0\smash \cp_0)}.
     \]
\end{cor}

\begin{proof}
Compare the calculation in Proposition \ref{prop:MUtheta_sigmap} with Proposition \ref{prop:calculate_K_KU}. (Note that the sign arises from the conventions used in defining the cobar complex, as in Proposition \ref{prop:cocycle-canonical-splitting}.)
\end{proof}

\section{Restricting to primitives}

The spherical elements of $MU_* (\cp_0\smash \cp_0)$ lie in the comodule primitives; this motivates the study of the restriction of the algebraic double transfer to the comodule primitives.

\subsection{Comodule primitives}

\begin{nota}
	For $M$ a left $MU_*MU$-comodule, write $\prim M$ for the graded abelian group of comodule primitives.
\end{nota}

As above $\hq$ denotes the rational Eilenberg-MacLane spectrum; its integral counterpart is denoted $\hz$. The following is clear:

\begin{lem}
\label{lem:thom_orientation_primitives}
	For $X$  a spectrum, there is a natural commutative diagram of graded abelian groups, induced by the orientation of $\hz$ and rationalization
\[
	\xymatrix{
\prim MU_* X
\ar@{^(->}[r]
\ar[d]
&
MU_* X
\ar[r]
&
\zed \otimes_{MU_* } MU_* (X)
\ar[r]
&
\hz_* X
\ar[d]
\\
\prim MU_* X \otimes \rat
\ar[rrr]_\cong
&&&
\hq_* X,
}
\]
in which the lower horizontal morphism is an isomorphism.

If $MU_* X$ has no additive torsion, then
$\prim MU_* X \rightarrow \hz_* X$ is a monomorphism.
\end{lem}

\begin{exam}
\label{exam:prim_CP}
	Let  $d$  be a natural number. Then
$$
	\prim
MU_* ((\cp_0)^{\smash d})
\hookrightarrow
\hz_* ((\cp_0)^{\smash d})
$$
is a morphism of algebras, where the product is induced by the $H$-space structure of $\cp$. The homology $\hz_*
((\cp_0)^{\smash d})$ is the free divided power algebra $\Gamma^* (\zed^{\oplus d}) $, hence $ \prim MU_* ((\cp_0)^{\smash d})$ is a subalgebra of $\Gamma^* (\zed^{\oplus d}) $.
\end{exam}

	The primitives $ \prim MU_* (\cp_0)$ were calculated by David Segal
\cite{segal}; an elegant approach is given by Miller in \cite[Proposition 4.1]{miller}, where
the primitive generators $p_n \in MU_{2n}(\cp_0)$ are defined by means of the expansion
\[
	\underline{\beta}_0(\exp (T))
=
\sum
\frac{p_n}{n!}T^n
\]
in $MU_* (\cp_0) \otimes \rat$, so that $|p_n|= 2n$.  The morphism
$
\prim MU_* (\cp_0) \rightarrow \hz_* (\cp_0)
$
sends $p_n$ to $ (\beta_1^\hz)^n  = n! \beta_n^\hz$, where  the $\beta_i^\hz$ denote the canonical module generators of $\hz_* (\cp_0)$ (cf. \cite[Remark 4.2]{miller}).

By the Künneth isomorphism $MU_* (\cp_0 \smash \cp_0) \cong MU_* (\cp_0) \otimes_{MU_*} MU_*(\cp_0)$, for pairs of natural numbers $(i,j)$, $p_i \otimes p_j$ is a primitive of $MU_* (\cp_0 \smash \cp_0)$.  The integral calculation of $\prim MU_* (\cp_0 \smash \cp_0)$   is an interesting and difficult problem: the elements $p_i \otimes p_j$ do not generate the primitives, due to delicate divisibility questions (cf. \cite{knapp_habilit}, \cite{baker_transfer} and \cite{BCRS}, for example).

\begin{lem}
\label{lem:identify-primitives}
	The primitive $\prim MU_* (\cp_0 \smash \cp_0)$ subgroup is a graded free
$\zed$-module such that $\prim MU_* (\cp_0 \smash \cp_0) \otimes \rat$ has
basis $\{ p_i \otimes p_j | i, j \geq 0 \}$.
\end{lem}

\begin{nota}
\label{nota:primel}
 Let $\primel (S, T) $ denote the two-variable power-series in $MU_* (\cp_0\smash \cp_0 ) \otimes \rat [[S, T]]$:
\[
\primel(S, T):= \sum_{m,n \geq 0} p_m \otimes p_n \frac{S^m T^n}{m! n !}.
\]
\end{nota}

\begin{lem}
 \label{lem:identify_primel}
There is an identity of formal power series:
\[
 \primel(S, T) = \underline{\beta}_0 (\exp S) \otimes \underline{\beta}_0 (\exp T).
\]
\end{lem}

\subsection{Restricting the double transfer to primitives}
\label{subsect:primitive_double_transfer}

 A primitive $\primel$ of degree $2k$ in $MU_* (\cp_0 \smash \cp_0)$ corresponds to
a morphism of comodules
\[
\primel \in
\hom_{MU_* MU}(MU_* [2k] , MU_* (\cp_0 \smash \cp_0)).
\]
This induces a class $\primel^* [\kappa] \in \ext^1_{MU_*MU} (MU_* [2k] , MU_* \otimes \qz [-4])$ and, by the chromatic connecting morphism $\partial_1$, the image of the double algebraic transfer:
\[
	\primel ^* [e_\tau^2] = \partial_1 \primel^*  [\kappa] \in \ext^2 _{MU_*MU} (MU_*[2k] , MU_* [-4]),
\]
where the identification follows from Proposition \ref{prop:kappa-transfer}.

In particular, to understand the restriction of the double algebraic transfer to the primitive element $\primel$, it suffices to consider $\primel^*  [\kappa]$,
which is represented by the cocycle $\kappa\circ \primel$:
\[
\xymatrix{
	MU_* [2k]
\ar[r]^(.4){\primel}
&
MU_*(\cp_0 \smash \cp_0 )
\ar[r]
^\kappa
&
MU_* MU \otimes \qz [-4].
}
\]
By Proposition \ref{prop:calculate_Xi},  the cocycle $\kappa\circ \primel$ fits
 into the commutative diagram
\[
\xymatrix{
	MU_* [2k]
\ar[r]^(.4){K \circ \primel}
\ar[rd]_{\kappa \circ \primel}
&
MU_*MU \otimes \rat [-4]
\ar@{->>}[d]
\\
&
MU_* MU \otimes \qz [-4].
}
\]
Recall that the morphism $K$ is given in Proposition \ref{prop:calculate_Xi} by specifying the image of $\underline{\beta}_0 (x) \otimes \underline{\beta}_0 (y)$.

\begin{nota}
	Write $\divBern_n ^L, \divBern_n^R \in MU_*MU \otimes \rat$ for the reduced Bernoulli numbers  associated to the left (respectively right) $MU_*$-algebra structures.
\end{nota}

\begin{prop}
 \label{prop:restrict_Xi_prim}
The restriction of the morphism $K : MU_* (\cp_0 \smash \cp_0) \rightarrow MU_*MU\otimes \rat [-4]$ to the primitive elements is determined
by
\[
 K (p_m \otimes p_n)
=
 (\divBern^R_{m+1} - \divBern^L _{m+1})\divBern^R_{n+1}
\]
for natural numbers $m,n$.
\end{prop}

\begin{proof}
By Proposition \ref{prop:calculate_Xi}, the morphism $K $ is given by
\[
 \underline{\beta}_0 (x) \otimes \underline{\beta}_0 (y)
\mapsto
\Big(
\frac{1}{x}- \frac{1}{\underline{b}(x)}
\Big)
\Big(
\frac{1}{\log^L y }
- \frac{1}{\underline{b}(y)}
\Big).
\]
 Lemma \ref{lem:identify_primel} identifies the generating formal power series  $\primel(S,T)$ for the primitive elements $p_i \otimes p_j$; thus the image of $\primel(S,T)$ is given by substituting the power series $x= \exp^L
S $, $y=\exp^L T$ in the above exression (note that  the left module structure of $MU_*MU$ is used), which gives
\[
\primel(S,T)
\mapsto
\Big(
\frac{1}{\exp^L S}- \frac{1}{\underline{b}(\exp^L S)}
\Big)
\Big(
\frac{1}{\log^L ( \exp^L T)}
- \frac{1}{\underline{b}(\exp^L T)}
\Big).
\]
Simplifying and reversing the order in the two brackets, this gives:
\[
\Big(
\frac{1}{\exp^R S}- \frac{1}{\exp^L S}
\Big)
\Big(
\frac{1}{\exp^R T}
-
\frac{1}{T}
\Big).
\]
The result follows from the definition of $\primel(S,T)$ and  of the reduced Bernoulli numbers.
\end{proof}

	In principle, this result determines the class $\primel ^* [\kappa]$, for any primitive $\primel$. This can be made more concrete by passing to elliptic homology and appealing to the invariants introduced by Laures \cite{laures} and by Behrens \cite{behrens}, as explained in the following sections.

\section{Passage to elliptic homology}

To study the $p$-local Adams-Novikov two-line, for $p \geq 5$ a prime, complex cobordism can usefully be replaced by elliptic homology, by change of rings.
 The results of this section are entirely algebraic, relying on the fact that the formal group law of elliptic homology is defined over the ring of holomorphic modular forms.

\subsection{Formal group law input}
\label{subsect:fgl_input}

Consider the ring $\mf$ of
holomorphic forms of level one over the ring $\zedsix$, so that $\mf \cong \zedsix [c_4, c_6]$,
 graded by weight. The Fourier expansion at the cusp at infinity defines the
$q$-expansion $\mf\hookrightarrow \zedsix [[q]] [u]$, a
monomorphism of rings by the $q$-expansion principle (the variable $u$
corresponds to the grading). The morphism $q^0 :
\zedsix[[q]] \rightarrow \zedsix$ sending a  power series to its
constant term induces a ring morphism $\zedsix [[q]] [u] \rightarrow \zedsix [u]
$. This gives a commutative diagram of ring morphisms:
\begin{eqnarray}
\label{eqn:q0_hol}
	\xymatrix{
L
\ar[r]
\ar[d]_{\gmu}
&
\mf
\ar[d]
\\
\zedsix [u]
&
\zedsix [[q]][u] .
\ar[l]^{q^0}
}
\end{eqnarray}
The composite $L \rightarrow \zedsix [[q]][u]$ classifies the  (graded) formal
group law  associated to the Tate Weierstrass curve defined
over $\zedsix [[q]] [u]$. After reduction $q \mapsto 0$, this is the
multiplicative (graded) formal group law $\gmu$ over $\zedsix [[q]][u]$, since the Tate curve has
multiplicative reduction.

\begin{rem}
	The universal Weierstrass curve over the ring of holomorphic modular forms $\mf$ is not an elliptic curve;
however, the associated formal group law is
defined, since the curve is smooth at the identity section. Similarly for  the Tate curve.
\end{rem}

The ring $\mfmero$ of  meromorphic modular forms over $\zedsix$ is isomorphic
to $ \mf [\Delta^{-1}]$, where $\Delta$ is the discriminant. The $q$-expansion of a
meromorphic modular form lies in $\zedsix [[q]][q^{-1}]$ and $q^0$ defines an
additive morphism $q^0 : \zedsix [[q]][q^{-1}] \rightarrow \zedsix$. This
gives the following commutative diagram
\begin{eqnarray}
\label{eqn:mf-fgl}
	\xymatrix{
L
\ar[r]
\ar[d]_{\gmu}
&
\mf
\ar@{^(->}[r]
\ar[d]
&
\mfmero
\ar[d]
\\
\zedsix [u^{\pm 1}]
&
\zedsix [[q]][u]
\ar@{^(->}[r]
\ar[l]^{q^0}
&
\zedsix [[q]][q^{-1},u^{\pm 1}],
\ar@{-->}@/^2pc/[ll]^{q^0}
}
\end{eqnarray}
in which the solid arrows are ring morphisms.

The morphisms $L \rightarrow \zedsix [u^{\pm 1}] $ and $L \rightarrow
\mfmero$ are Landweber exact and correspond respectively to
 $KU_\zedsix$ (which will be denoted here simply by $KU$) and a version of elliptic homology, denoted  by  $\tmf$.

\begin{prop}\cite[Theorem 2.7]{laures}.
The ring $\tmfcoeff _ 0 KU$ is isomorphic to Katz's  ring
$\divcong_\zedsix$ of divided congruences over $\zedsix$.
\end{prop}

\begin{rem}
\
\begin{enumerate}
\item
The ring $\tmfcoeff_* KU $ is concentrated in even degree and is $2$-periodic.
	\item
The ring $\tmfcoeff_0 KU$ is a subring of $\tmfcoeff_0 KU \otimes \rat\cong
(\tmfcoeff_* \otimes KU_* \otimes \rat)_0$ and  identifies with the sub-$\zedsix$-module
of sums $\sum_i f_i $ of modular forms such that the Fourier expansion $\Sigma_i f_i(q)$ has coefficients
in $\zedsix$; this is precisely the ring $\divcong_\zedsix$ of divided congruences \cite{katz}.
\end{enumerate}
\end{rem}

\subsection{The reduction map $\rhobar^1$}
\label{subsect:rho}

The additive morphism $q^0 : \mfmero \rightarrow \zedsix [u^{\pm 1}]$ is
realized by a morphism of spectra $q^0 : \tmf \rightarrow KU$, which is
derived from Miller's elliptic character (see
\cite{laures}, following Miller \cite{millerW}). Hence there is an induced
morphism of spectra
$
	\tmf \smash \tmf
\rightarrow
KU \smash \tmf ,
$
which induces a morphism of right $\tmfcoeff_*$-modules
\[
	\rhobar^1 : \tmfcoeff_* \tmfcoeff
\rightarrow
KU_* \tmfcoeff,
\]
which is used in defining Laures' $f$-invariant \cite{laures}  (see Section \ref{subsect:f_inv}).

Since  $\tmf$ and $KU$ are Landweber exact,  $\tmfcoeff_* \tmfcoeff \otimes \rat \cong
\tmfcoeff_* \otimes \tmfcoeff_* \otimes \rat$ and $KU_* \tmfcoeff \otimes
\rat \cong  KU_*\otimes \tmfcoeff_*  \otimes \rat$.

\begin{prop}
\label{prop:rho_rat}
The morphism $\rhobar^1 \otimes \rat$ is the morphism of right $\tmfcoeff_* \otimes \rat$-modules:
\[
q^0 \otimes \tmfcoeff_*  \otimes \rat :
\tmfcoeff_* \otimes \tmfcoeff_* \otimes \rat
\cong
\tmfcoeff_* \tmfcoeff \otimes \rat
\rightarrow
 KU_* \tmfcoeff\otimes \rat
 \cong
KU_* \otimes  \tmfcoeff_* \otimes \rat
\]
\end{prop}

There is a morphism of short exact sequences of right  $\tmfcoeff_*$-modules
\begin{eqnarray}
\label{eqn:ses-rho1}
	\xymatrix{
\ \ \ \ 0
\ar[r]
&
\tmfcoeff_*\tmfcoeff
\ar[r]
\ar[d]_{\rhobar^1}
&
\tmfcoeff_*\tmfcoeff \otimes \rat
\ar[r]
\ar[d]^{\rhobar^1 \otimes \rat}
&
\tmfcoeff_*\tmfcoeff \otimes \qzsix
\ar[r]
\ar[d]^{\rhobar^1 \otimes \qzsix}
&
0
\\
0
\ar[r]
&
KU_* \tmfcoeff
\ar[r]
&
KU_* \tmfcoeff\otimes \rat
\ar[r]
&
KU_* \tmfcoeff \otimes \qzsix
\ar[r]
&
0.
}
\end{eqnarray}
Thus  $\rhobar^1 \otimes \qzsix$ is determined by $\rhobar^1 \otimes \rat$ and
hence by the additive morphism $q^0$.

\begin{rem}
	Following \cite[Theorem 4.2]{behrens_laures}, the morphism $\rhobar^1$ is used here to define the $f$-invariant rather than the analogous morphism
	$\rho^1 : \tmfcoeff_*\tmfcoeff \rightarrow \tmfcoeff_* KU$ (compare  \cite[Proposition 3.9]{laures}). The relationship between the two approaches to calculating the $f$-invariant is explained by \cite[Proposition 3.10]{laures}.
\end{rem}

\subsection{Reduction of the cocycle  $\kappa$}
\label{subsect:basechange_kappa_ell}

The class $[\kappa] \in \ext^1 _{MU_*MU} (MU_* (\cp_0 \smash \cp_0), MU_* \otimes \qz[-4])$ corresponds to the algebraic double transfer, by Proposition \ref{prop:kappa-transfer} and, by  Proposition \ref{prop:calculate_Xi}, the representing  cocycle
\[
\kappa \in \hom_{MU_*} (MU_* (\cp_0 \smash \cp_0) , MU_*MU \otimes \qz[-4])
\]
 is induced by the morphism $K \in \hom_{MU_*} (MU_*
(\cp_0 \smash
\cp_0) , MU_*MU \otimes \rat [-4])$, by composition with the quotient map $MU_*MU
\otimes \rat \twoheadrightarrow MU_*MU \otimes \qz$.

Base change along $MU_* \rightarrow \tmfcoeff_*$ gives a cocycle
\[
\kappa_{\tmf} :=\  _{\tmfcoeff_*}\kappa _{\tmfcoeff_*} \in \hom_{\tmfcoeff_*}
(\tmfcoeff_*(\cp_0 \smash \cp_0) ,
\tmfcoeff_* \tmfcoeff \otimes \qzsix [-4]),
\]
which represents a class $[\kappa_{\tmf}] \in \ext^1 _{\tmfcoeff_*\tmfcoeff} (\tmfcoeff_* (\cp_0 \smash \cp_0), \tmfcoeff_* \otimes \qz[-4])$,
by Lemma \ref{lem:cocycle-base-change}. The following is clear:

\begin{lem}
\label{lem:kappa_ell_Xi}
The morphism $\kappa_{\tmf}$ is the reduction of the morphism
\[
_{\tmfcoeff_*}K_{\tmfcoeff_*} \in
 \hom_{\tmfcoeff_*} (\tmfcoeff_*(\cp_0 \smash \cp_0) ,
\tmfcoeff_* \tmfcoeff \otimes \rat [-4])
\]
 via the morphism
$\tmfcoeff_* \tmfcoeff \otimes \rat [-4] \twoheadrightarrow \tmfcoeff_* \tmfcoeff\otimes \qzsix
[-4]$.
\end{lem}

\begin{prop}
\label{prop:kappa-tmf}
	The morphism of right $\tmfcoeff_*$-modules
\[
(\rhobar^1 \otimes \qzsix ) \circ \kappa_\tmf \in \hom
(\tmfcoeff_*(\cp_0 \smash \cp_0) ,
KU_* \tmfcoeff \otimes \qzsix [-4])
\]
 coincides with the morphism $_{KU_*} \kappa _{\tmfcoeff_*}$.

Hence, the morphism $(\rhobar^1 \otimes \qzsix ) \circ \kappa_\tmf$ is the reduction
of the morphism of right $\tmfcoeff_*$-modules
\[_{KU_*}K_{\tmfcoeff_*} \in \hom (\tmfcoeff_*
 (\cp_0\smash \cp_0) , KU_* \tmfcoeff\otimes \rat [-4])
\]
 via the morphism
$KU_* \tmfcoeff\otimes \rat [-4] \twoheadrightarrow KU_* \tmfcoeff \otimes \qzsix
[-4]$.
\end{prop}

\begin{proof}
	The diagram of short exact sequences (\ref{eqn:ses-rho1}) together
with Lemma \ref{lem:kappa_ell_Xi} show that it is sufficient to calculate the
respective morphisms to $KU_* \tmfcoeff \otimes \rat$. This can be carried out
using the identification of $\rhobar^1 \otimes \rat$ given by Proposition \ref{prop:rho_rat}.

There is a commutative diagram
\[
 \xymatrix{
\tmfcoeff_* (\cp_0 \smash \cp_0)
\ar@/^1pc/[rd]^{_{\tmfcoeff_*}K_{\tmfcoeff_*}}
\ar[d]_{_{\mf}K_{\tmfcoeff_*}}
\\
\mf \otimes_{MU_*} MU_*\tmfcoeff  \otimes \rat [-4]
\ar[d]_{q^0 \otimes \tmfcoeff_*}
\ar[r]
&
\tmfcoeff_* \otimes _{MU_*} MU_* \tmfcoeff_*\otimes
\rat [-4]
\ar@/^1pc/[ld]^(.4){\rhobar^1 \otimes \rat}
\\
KU_* \tmfcoeff \otimes \rat [-4],
}
\]
where the horizontal morphism is induced by $\mf \rightarrow  \mfmero \cong \tmfcoeff_*$.
The commutativity of the top triangle follows from the fact that the
elliptic formal  group law  is defined over the ring of holomorphic modular forms, $\mf$, and the commutativity of
the lower triangle follows from the commutative diagram (\ref{eqn:mf-fgl}).

To complete the proof, observe that the vertical composite is the morphism $_{KU_*}K_{\tmfcoeff_*}$, by the
commutative diagram (\ref{eqn:q0_hol}).
\end{proof}

\section{The $f$ and $f'$ invariants}
\label{sect:invariants}

The algebraic image of the double transfer can be analysed by using either the $f$-invariant of Laures (considered as an
invariant of the Adams-Novikov two-line) or the $f'$-invariant introduced by Behrens.

\subsection{Recollections on the $f$ invariant}
\label{subsect:f_inv}

The $f$-invariant of Laures  \cite{laures} is a
homomorphism
\[
	f :
\pi_{2k}(S) \otimes \zedsix
\rightarrow
\divcong _\rat / \Big( \divcong_\zedsix + (\mfmero_0)_\rat +
(\mfmero_{k+1})_\rat \Big).
\]
This factorizes across an invariant
\[
\iota^2 : 	\ext^{2,2k+2} _{MU_*MU} (MU_* , MU_*) \otimes \zedsix
\hookrightarrow
\divcong _\rat / \Big(\divcong_\zedsix \oplus (\mfmero_0)_\rat \oplus
(\mfmero_{k+1})_\rat \Big),
\]
where the injectivity is given by \cite[Proposition 3.9]{laures}.

Via  the chromatic connecting map
\[
	\ext^{1,*}_{MU_*MU} (MU_*,MU_* \otimes \qzsix )
\stackrel{\partial_1}{\rightarrow}
 \ext^{2,*}_{MU_*MU} (MU_*,MU_* ) \otimes \zedsix,
\]
$\iota^2$ defines an invariant of
$
\ext^{1,*}_{MU_*MU} (MU_*,MU_* \otimes \qzsix ).
$

Change of rings associated to the orientation $MU_* \rightarrow \tmfcoeff_*$, allows the respective groups to be  replaced  by
\begin{eqnarray*}
&&	\ext^{2,*} _{\tmfcoeff_*\tmfcoeff} (\tmfcoeff_* , \tmfcoeff_*)
	\\
&&	\ext^{1,*}_{\tmfcoeff_*\tmfcoeff} (\tmfcoeff_*,\tmfcoeff_* \otimes \qzsix ).
\end{eqnarray*}

We identify the invariant $\iota^2$ following Behrens and Laures. Write $M_k ^{\bullet +1} \cong \pi_{2k} (\tmf^{\bullet +1})$ for the cobar complex associated to
$\tmf$. A morphism between semi-cosimplical abelian groups is defined \cite[page 25]{behrens_laures}:
\begin{eqnarray}
\label{eqn:BL_semicosimplicial}
	\xymatrix{
M_k^{(1)}
\ar@<.5ex>[r]|{d_0}
\ar@<-.5ex>[r]|{d_1}
\ar[d]
_{\overline{\rho}^0}
&
M_k^{(2)}
\ar@<1ex>[rr]|{d_0}
\ar@<0ex>[rr]|{d_1}
\ar@<-1ex>[rr]|{d_2}
\ar[d]
_{\overline{\rho}^1}
&&
M_k^{(3)}
\ar[d]
\ar[d]
_{\overline{\rho}^2}
\\
(\mfmero_k)_{\zedsix}
\ar@<.5ex>[r]|(.6){d_0}
\ar@<-.5ex>[r]|(.6){d_1}
&
\divcong_{\zedsix}
\ar@<1ex>[rr]|(.3){d_0}
\ar@<0ex>[rr]|(.3){d_1}
\ar@<-1ex>[rr]|(.3){d_2}
&&
\divcong_\rat / \big(
(\mfmero_k)_{\rat}
\oplus
(\mfmero_0)_\rat
\big)
}
\end{eqnarray}

The composite of $\overline{\rho}^2$ with the projection
\[
	\divcong_\rat / \big(
(\mfmero_k)_{\rat}
\oplus
(\mfmero_0)_\rat
\big)
\twoheadrightarrow
\divcong_\rat / \big(
\divcong_{\zedsix} \oplus (\mfmero_k)_{\rat} \oplus (\mfmero_0)_\rat
\big)
\]
induces the morphism
\[
\iota^2 : 	\ext^{2,2k} _{\tmfcoeff_*\tmfcoeff} (\tmfcoeff_* , \tmfcoeff_*) \otimes \zedsix
\hookrightarrow
\divcong _\rat / \Big(\divcong_\zedsix \oplus (\mfmero_0)_\rat \oplus
(\mfmero_{k})_\rat \Big),
\]
 on restriction to cocycles.

Write the  chain cocomplex associated to the cobar complex as
\[
	\xymatrix{
M_k^{(1)}
\ar[r]^{\delta^0}
&
M_k^{(2)}
\ar[r]^{\delta^1}
&
M_k^{(3)}
\ar[r]
&
\ldots
\ \  .
}
\]
The morphism $\iota^2$ is identified explicitly by the following straightforward application of chromatic arguments.

\begin{lem}
\label{lem:iota2_rho1}
	Let $x$ be a $2$-cocycle in $M_k^{(2)}$ which represents a class $$[x] \in \ext^{2, *} _{\tmfcoeff_* \tmfcoeff} (\tmfcoeff_*, \tmfcoeff_*).$$ Then
	\begin{enumerate}
		\item
		there exists an element $c \in M_k^{(1)}$ and an integer $n$ such that $\delta^1 c = nx$;
		\item
		the invariant $\iota^2 [x]$ is represented by the element $\frac{1}{n}\overline{\rho}^1 (c) \in \divcong_\rat$.
	\end{enumerate}
\end{lem}

\begin{proof}
	A straightforward consequence of the commutative diagram (\ref{eqn:BL_semicosimplicial}) together with the fact that $\ext^d_{\tmfcoeff_*\tmfcoeff}(\tmfcoeff_* , \tmfcoeff_*) \otimes \rat$ is trivial
for $d >0$.
\end{proof}

\begin{prop}
\label{prop:f-inv_1cocycle}
	Let $[c] \in \ext^{1,2k}_{\tmfcoeff_*\tmfcoeff} (\tmfcoeff_* , \tmfcoeff_* \otimes \qzsix)$ be represented by a
cocycle $c : \tmfcoeff_* [2k] \rightarrow \tmfcoeff_* \tmfcoeff \otimes \qzsix$ which factorizes as
left $\tmfcoeff_*$-module morphisms
\[
	\xymatrix{
\tmfcoeff_*[2k]
\ar[r]^(.4){\hat{c}}
\ar[rd]
_c
&
\tmfcoeff_*\tmfcoeff \otimes \rat
\ar@{->>}[d]
\\
&
\tmfcoeff_*\tmfcoeff \otimes \qzsix.
}
\]
Then the invariant  $\iota^2( \partial_1 [c])$ is represented by the image of the generator
under the map
\[
\xymatrix{
	\tmfcoeff_*[2k]
\ar[rr]^(.45){\ _{KU_*} \hat{c} _{\tmfcoeff_*} }
&&
KU_*\tmfcoeff \otimes \rat,
}
\]
where $KU_{2k} \tmfcoeff\otimes \rat$ is identified with $\divcong_\rat$  by periodicity.
\end{prop}

\begin{proof}
The result follows from Lemma \ref{lem:iota2_rho1}.
\end{proof}

\subsection{Restricting the $f$-invariant to primitives}

\begin{nota}
For $\primel$ a primitive in $\prim MU_* (\cp_0 \smash \cp_0)$, let $f(\primel)$
denote the $f$-invariant of $ \primel^* [\kappa]\in \ext^1_{MU_*MU}
(MU_*[|\primel| +4] , MU_* \otimes \qzsix)$.
\end{nota}

For $n$ a natural number, write $\divBern_n ^{\tmf}\in \tmfcoeff_* \otimes \rat$ (respectively
$\divBern_n ^{KU} \in KU_* \otimes \rat$) for the reduced Bernoulli numbers associated to the complex
orientations of $\tmf$ and $KU$ respectively.

\begin{rem}
	The reduced Bernoulli number $\divBern_n ^{\tmf}$ is defined in the ring $\mf \otimes \rat$ of holomorphic modular forms, since the formal group law of $\tmfcoeff_*$ is
the image of a formal group law over $\mf$ via the morphism $\mf\hookrightarrow \mfmero\cong \tmfcoeff_*$.
\end{rem}

\begin{thm}
\label{thm:f-pspt}
Let  $s,t$ be natural numbers.
The $f$-invariant $$f(p_s \otimes p_t)
\in \divcong _\rat / \Big(\divcong_\zedsix \oplus (\mfmero_0)_\rat \oplus
(\mfmero_{s+t+2})_\rat \Big),$$  is represented by the element
$
- \divBern_{t+1}^\tmf \divBern _{s+1} ^{KU}
\in
\divcong_\rat.
$
\end{thm}

\begin{proof}
By Proposition \ref{prop:f-inv_1cocycle}, these invariants are represented by the morphism
$
	\ _{KU_*} ( K \circ \primel) _{\tmfcoeff_*}.
$
Hence, by  Proposition \ref{prop:restrict_Xi_prim},
 $f(p_s \otimes p_t)$ is represented by
\[
	 (\divBern_{s+1} ^{\tmf} - \divBern _{s+1} ^{KU}) \divBern_{t+1}^{\tmf}
	\in
\divcong_\rat.
\]
The term $\divBern_{s+1}^{\tmfcoeff} \divBern_{t+1} ^{\tmfcoeff}$ becomes zero on passage to the quotient, since it belongs to
the subgroup $(\mfmero_{s+t+2})_\rat$.
\end{proof}

\begin{cor}
\label{cor:symm_f-invariant}
Let  $s,t$ be natural numbers.
The invariant $f(p_s \otimes p_t)$  is represented by the element
$
 \divBern_{s+1}^\tmf \divBern _{t+1} ^{KU}
\in
\divcong_\rat.
$
\end{cor}

\begin{proof}
The group $\symm_2$ acts on $MU_* (\cp_0 \smash \cp_0)$ by comodule morphisms induced by interchanging the factors $\cp_0$. It is straightforward to show that the induced right action on $\ext_{MU_*MU} ^2 (MU_* (\cp_0 \smash \cp_0) , MU_* [-4]) $ satisfies
\[
	([e_\tau]^2 ) \sigma = \mathrm{sgn} (\sigma) [e_\tau]^2.
\]
Hence
$
	f (p_s \otimes p_t) = - f (p_t \otimes p_s);
$
 in particular, the $f$-invariant of $p_s \otimes p_t$ is  represented by the element $\divBern_{s+1}^\tmf \divBern _{t+1} ^{KU}$.
\end{proof}

\subsection{The $f'$-invariant on primitives}

For $p \geq 5$ a prime, Behrens \cite{behrens} defines the $f'$-invariant via a morphism
\[
	f' :
\ext^{2,2k+2} _{MU_*MU} (MU_* , MU_*)_{(p)}
\rightarrow
H^0(C (l)^\bullet/p^\infty, v_1 ^\infty )_{2k+2},
\]
where $l$ is a topological generator of $\padics^\times$ (for example $l = \gamma$) and $C(l)^\bullet$ is an
explicit semi-cosimplicial abelian group which is defined in terms of modular forms of level one and modular forms of level $l$. Namely, as in \cite{behrens}, write $M_k (\Gamma_0 (l))_{\padics}$ for the space of modular forms of weight $k$ and level $\Gamma_0(l)$  over $\padics$
which are meromorphic at the cusps. Then the semi-cosimplicial graded abelian group is of the form
\[
C(l)^\bullet_{2k} =
\Big(
\xymatrix{
(M_k)_{\padics}
\ar@<.5ex>[r]|(.4){d_0}
\ar@<-.5ex>[r]|(.4){d_1}
&
*{
\begin{array}{c}
M_k(\Gamma_0 (l))_{\padics}\\
\times\\
	(M_k)_{\padics}
\end{array}
}
\ar@<1ex>[r]
\ar[r]
\ar@<-1ex>[r]
&
M_k(\Gamma_0 (l))_{\padics}
	}
\Big),
\]
where the morphisms $d_0, d_1$ are identified explicitly in terms of $q$-expansions.
(See \cite[Section 6]{behrens} and the review in \cite[Section 3]{behrens_laures}.) It follows that the $f'$-invariant of a class is
represented by a modular form which satisfies certain congruences.

\begin{rem}
	Behrens and Laures work $p$-locally and replace $MU$ by $BP$, so as to accord better with the results of Miller, Ravenel and Wilson \cite{MRW}. Thus, below  $\tmfcoeff_*$
 denotes $p$-local elliptic homology ($p \geq 5$), and a $p$-typical orientation $BP_* \rightarrow \tmfcoeff_*$ is fixed, as in \cite{behrens_laures}.
\end{rem}

Behrens and Laures shows that the $f'$-invariant fits into a commutative diagram
\begin{eqnarray}
	\label{eqn:f'_invariant}
	\xymatrix{
\ext^{2,4t}_{BP_*BP} (BP_*, BP_*)
\ar[dd]_{f'}
&
\ext^{0,4t}_{BP_*BP} (BP_*, BP_*/p^\infty,v_1^\infty)
\ar[l]_{\partial_1\partial_2}
\ar[d]^{L_{v_2}}
\ar@/_1pc/[ldd]
\\
&\ext^{0,4t}_{BP_*BP} (BP_*, BP_*/p^\infty,v_1^\infty [v_2^{-1}])
\ar[d]^\cong
\\
H^0(C (l)^\bullet/p^\infty, v_1 ^\infty )_{4t}
&
\ext^{0,4t}_{\tmfcoeff_*\tmfcoeff} (\tmfcoeff_*, \tmfcoeff_*/p^\infty,v_1^\infty),
\ar[l]^(.55){\cong}_(.55){\tilde{\eta}}
	}
\end{eqnarray}
where the diagonal arrow is induced  by the $p$-typical orientation of $\tmf$,  the vertical change of rings isomorphism is given in the proof of
\cite[Lemma 4.6]{behrens_laures} and the  isomorphism $\tilde{\eta}$ in \cite[Proposition 3.17]{behrens_laures}.  The upper triangle is commutative by
\cite[Diagram 3.15]{behrens_laures} and the lower triangle is commutative by \cite[Diagram 3.16]{behrens_laures}. Up to the isomorphism $\tilde{\eta}$, the $f'$-invariant can be considered
as taking values in the comodule primitives of $\tmfcoeff_*/p^\infty,v_1^\infty$; Behrens gives a modular description of $H^0(C (l)^\bullet/p^\infty, v_1 ^\infty )_{4t}$
 in  \cite[Theorems 1.2, 1.3]{behrens}.

Consider classes which are in the image of the algebraic double transfer. Recall that $\overline{\Theta}$ defines a chromatic factorization of the double transfer and this
induces a commutative diagram
\begin{eqnarray}
	\label{eqn:alg_double_transfer}
\xymatrix{
&
\ext_{BP_*BP}^{0,*-4} (BP_* , BP_* (\cp_0\smash \cp_0))
\ar[d]^{BP_* (\overline{\Theta})}
\ar[ld]
\\
\ext^{2,*}_{BP_* BP} (BP_* , BP_*)
&
\ext^{0,*}_{BP_* BP} (BP_* , BP_* /p^\infty, v_1^\infty)
\ar[l]^{\partial_1 \partial_2},
}
\end{eqnarray}
by Proposition \ref{prop:kappa-transfer} and Proposition \ref{prop:second-factorization} (with  $BP$ in place of $MU$).

\begin{lem}
\label{lem:independence_Theta}
	The composite morphism
	\[
		\xymatrix{
\ext_{BP_*BP}^{0,4t -4} (BP_* , BP_* (\cp_0\smash \cp_0) )
\ar[rr]^{BP_* (\overline{\Theta})}
&&
\ext^{0,4t}_{BP_* BP} (BP_* , BP_* /p^\infty, v_1^\infty)
\ar[d]
\\
&&
\ext^{0,4t}_{\tmfcoeff_*\tmfcoeff} (\tmfcoeff_*, \tmfcoeff_*/p^\infty,v_1^\infty)
}
	\]
induced by the change of rings associated to the $p$-typical orientation $BP_* \rightarrow \tmfcoeff_*$, is independent of the choice of $\overline{\Theta}$.
\end{lem}

\begin{proof}
	Follows from the commutativity of the diagrams  (\ref{eqn:f'_invariant}) and (\ref{eqn:alg_double_transfer}), together with the fact that the bottom horizontal morphism in  (\ref{eqn:f'_invariant})  is an isomorphism.
\end{proof}

By naturality, it suffices to replace the composite morphism considered above by the morphism
\[
\tmfcoeff_* (\overline{\Theta}) :
\ext_{\tmfcoeff_*\tmfcoeff}^{0,4t -4} (\tmfcoeff_* , \tmfcoeff_* (\cp_0\smash \cp_0) )
\rightarrow
\ext_{\tmfcoeff_*\tmfcoeff}^{0,4t} (\tmfcoeff_* , \tmfcoeff_* /p^\infty,v_1^\infty )
\]
which is considered as being the $f'$ invariant for primitive elements. This morphism is independent of the choice of orientation; hence, in the following, $\tmfcoeff_*$ is equipped with its standard complex orientation.

Theorem \ref{thm:MU-thetas} gives a commutative
diagram
\[
\xymatrix{
&
\prim \tmfcoeff_*(\cp_0 \smash \cp_0)
\ar@{^(.>}[d]
\ar@/^1pc/@<4em>@{.>}[dd]^{f'}
\ar@/^3pc/@<4em>@{.>}[ddd]^{f''}
\\
\tmfcoeff_* (\cp_{-1} \smash \cp_0)
\ar[dd]_{\tmfcoeff_* (\theta)}
&
\tmfcoeff_* (\cp_0 \smash \cp_0)
\ar@{.>}[l]_{\sigma'}
\ar[d]^{\tmfcoeff_*(\overline{\Theta})}
\\
&
\tmfcoeff_* /p^\infty, v_1 ^\infty [-4]
\ar@{^(->}[d]
\\
\tmfcoeff_*KU \otimes \rat[-4]
\ar[r]
&
\tmfcoeff_* KU \otimes \rat/\Big(\tmfcoeff_*KU _{(p)} \oplus (\tmfcoeff_*)_\rat \Big) [-4],
}
\]
where the solid arrows denote comodule morphisms and the dotted arrows morphisms of
$\plocal$-modules. The composite of $f'$ with the monomorphism $\tmfcoeff_* /p^\infty, v_1 ^\infty [-4]
\hookrightarrow \tmfcoeff_* KU \otimes \rat/\Big(\tmfcoeff_*KU _{(p)} \oplus (\tmfcoeff_*)_\rat \Big) [-4]$ is denoted $f''$, as indicated.

\begin{rem}
	The morphism $\tmfcoeff_* (\theta)$ composed with the morphism induced by $\psi^\gamma - 1$ is integral in the appropriate sense, as a consequence of the construction of $\theta$. This can be related to the analysis of $H^0 (C(l)^\bullet/ p^\infty, v_1^\infty) $ using the identifications of \cite[Proposition 6.1]{behrens}.
\end{rem}

The morphism $\tmfcoeff_*(\theta) \circ \sigma'$ is given by Proposition
\ref{prop:MUtheta_sigmap}, after base change  to $\tmfcoeff_*$; it is determined by $
 \underline{\beta}_0 (x) \otimes \underline{\beta}_0 (y)
\mapsto $
\[
\Big( \frac{1}{\underline{b}'(x) } - \frac{1}{x}\Big)
\Big(\frac{1}{\log^{\tmf} y} - \frac{1}{\underline{b}' (y) }\Big)
+
\sum _{i, j >0}\frac{B_i^{KU}}{i!}\frac{B_j^{KU}}{j!}
\left(
\frac{\gamma^i -1}{\gamma^{i+j}-1}
\right)
(\log^\tmf x)^{i-1}(\log ^\tmf y) ^{j-1},
\]
where the power series $\underline{b}'$ is understood as $\exp^{KU} \circ
\log^\tmf$ when considered as a power series over $\tmfcoeff_* KU  \otimes \rat$.

\begin{thm}
\label{thm:fprime}
	The $f'$-invariant on the primitives of $\tmfcoeff_* (\cp_0\smash \cp_0) $ is determined by
\[
f'' (p_s \otimes p_t)
 =
\Big[
 \divBern^\tmf_{s+1}\divBern^{KU} _{t+1}
  + \divBern^{KU}_{s+1}\divBern^{KU}
_{t+1}\frac{\gamma^{s+1}(1- \gamma^{t+1}) }{\gamma^{s+t+2}-1}
\Big],
\]
where $s, t$ are natural numbers and $
 \divBern^\tmf_{s+1}\divBern^{KU} _{t+1}  + \divBern^{KU}_{s+1}\divBern^{KU}
_{t+1}\frac{\gamma^{s+1}(1- \gamma^{t+1}) }{\gamma^{s+t+2}-1}
$ is considered as an element of $\tmfcoeff_*KU \otimes \rat$.
\end{thm}

\begin{proof}
The method of proof is similar to that used in Proposition \ref{prop:restrict_Xi_prim}.
Substitute $x= \exp^\tmf S$ and $y = \exp^\tmf T$ in the power
series representing $\tmfcoeff_* (\theta) \circ \sigma'$; this gives the power series
 \[
\Big(  \frac{1}{\exp^{KU} S} - \frac{1}{\exp^\tmf S  }\Big)
\Big(\frac{1}{T} - \frac{1}{\exp^{KU} T }\Big)
+
\sum _{i, j >0}\frac{B_i^{KU}}{i!}\frac{B_j^{KU}}{j!}
\left(
\frac{\gamma^i -1}{\gamma^{i+j}-1}
\right)
S^{i-1}T  ^{j-1}.
\]
The result follows by identifying coefficients.
\end{proof}

\begin{rem}
\label{rem:relate_f,f'}
	The relationship between the $f$ and the $f'$ invariants (in the general case) is made explicit in \cite[Theorem 4.2]{behrens_laures} by analysing the semi-cosimplicial diagram (\ref{eqn:BL_semicosimplicial}).

Upon restricting to classes arising from comodule primitives via the algebraic double transfer, the relationship is clear. Observe that
Theorem \ref{thm:f-pspt}, Corollary \ref{cor:symm_f-invariant} and Theorem \ref{thm:fprime} show that,  on passage to the quotient $$\divcong _\rat / \Big(\divcong_\zedsix \oplus (\mfmero_0)_\rat \oplus
(\mfmero_{s+t+2})_\rat \Big)$$ the elements $f(p_s \otimes p_t)$ and $f'(p_s \otimes p_t)$ both are defined by the class of the  element $\divBern^\tmf_{s+1}\divBern^{KU} _{t+1}  $ in $\divcong_\rat$, since the additional term appearing in Theorem \ref{thm:fprime} becomes trivial in this quotient.

The relationship can be seen as  a direct consequence of  Lemma \ref{lem:restrict_Theta} and Corollary \ref{cor:relate_MUtheta}; here, the sign appearing in the leading term in Corollary \ref{cor:relate_MUtheta}  has been avoided by following Behrens and Lawson in defining the $f$-invariant by using the morphism $\rhobar^1$.

\end{rem}

\appendix
\section{Cohomology and cocycles for Hopf algebroids}
\label{sect:cohom}

The reader is referred to  \cite[Appendix A2]{ravenel_green} for basics on Hopf algebroids.
For $(A,\Gamma)$  a Hopf algebroid, the category of left $\Gamma$-comodules is denoted here $\Gamma \dash \comod$ and the category of $A$-modules  $A\dash \modules$; the structure morphism of a $\Gamma$-comodule $M$ is written $\psi_M : M \rightarrow \Gamma \otimes_A M$.

\subsection{Recollections on cocycles}
\label{subsect:cocycles}

Throughout this section, $(A,\Gamma)$ is a flat Hopf algebroid, so that $\Gamma \dash \comod$ is abelian with enough injectives; extension groups in this category are denoted by $\ext_\Gamma$;  the extended comodule functor $\Gamma \otimes _A - : A\dash \modules \rightarrow \Gamma \dash \comod$ is right adjoint to the forgetful functor $\Gamma \dash \comod \rightarrow A \dash \modules$.

For $M$  a $\Gamma$-comodule, let $(\cobar_\Gamma ^\bullet M, d^\bullet)$
denote  the unreduced cobar resolution of $M$ (Cf. \cite[Definition A1.2.11]{ravenel_green}).

\begin{defn}
\
\begin{enumerate}
	\item
A $\Gamma$-comodule $X$ is $A$-projective if its underlying $A$-module is projective.
\item
	A $\Gamma$-comodule $J$ is relatively injective if it has the morphism extension problem
with respect to monomorphisms of $\Gamma$-comodules which are split as $A$-module morphisms.
\end{enumerate}
\end{defn}

\begin{lem}
	A $\Gamma$-comodule $J$ is relatively injective if and only if it is a direct summand of an extended comodule.
\end{lem}

\begin{lem}
\label{lem:ext_vanishing}
	Let $M$ be a $\Gamma$-comodule and $N$ an $A$-module. Then there is a natural isomorphism
\[
	\ext^i _\Gamma (M, \Gamma \otimes _A N )
\cong
\ext^i _A (M, N) .
\]
Hence, if $X$ is an $A$-projective $\Gamma$-comodule and $J$ is a relatively injective $\Gamma$-comodule,
\[
	\ext^i _\Gamma(X, J) = \left\{ \begin{array}{ll}
	                               	\hom_\Gamma (X, J) & i= 0
\\
0
& i >0.
	                               \end{array}
\right.
\]
\end{lem}

\begin{proof}
	The result follows from the argument of \cite[Lemma A1.2.8b)]{ravenel_green}.
\end{proof}

\begin{prop}
\label{prop:ext-cobar}
\cite[Lemma A.1.1.6, Corollary A1.2.12]{ravenel_green}
Let $X$ be an $A$-projective  $\Gamma$-comodule and $M$ be a  $\Gamma$-comodule. Then $ \ext^\bullet_\Gamma (X, M)$ can be calculated by a resolution of $M$ of the from $M \rightarrow J^\bullet$, where each $J^k$ is relatively injective.
 In particular, $ \ext^\bullet_\Gamma (X, M)$ is naturally isomorphic to the cohomology
of the cocomplex $\hom_\Gamma (X, \cobar^\bullet_\Gamma M)$, which has terms
$\hom_A (X, \Gamma^{\otimes_A s }\otimes _A M)$.
\end{prop}

\begin{cor}
\label{cor:1-cocycle}
Let  $X$ be an $A$-projective $\Gamma$-comodule and $M$ be a  $\Gamma$-comodule. An
extension class $[\kappa] \in \ext^1 _\Gamma (X, M)$ is represented by a cocycle
$\kappa \in
\hom_A (X, \Gamma \otimes_A M )$.
\end{cor}

The following statement makes explicit the first differential $d^0_M$ of the  complex $\hom_\Gamma (X, \cobar^\bullet_\Gamma M)$.

\begin{lem}
\label{lem:d0}
Let $X, M$ be as in Proposition \ref{prop:ext-cobar} and let $\alpha : M \rightarrow N$ be a morphism of $\Gamma$-comodules.
\begin{enumerate}
	\item
The differential
$d^0_M : \hom_A (X, M) \rightarrow \hom_A (X, \Gamma \otimes_A M )$ is given by
 $
d^0_M g := (\Gamma \otimes g ) \circ \psi_X - \psi_M \circ g.
$
	\item
	The differential
$d^0_N : \hom_A (X, N) \rightarrow \hom_A (X, \Gamma \otimes_A N )$ satisfies the relation
$
d^0_N (\alpha \circ g) = (\Gamma \otimes _A \alpha) \circ d^0_M g.
$
\end{enumerate}
\end{lem}

\begin{lem}
\label{lem:explicit-connecting}
	Let  $X$ be an $A$-projective $\Gamma$-comodule and
\[
	0
\rightarrow
N_1
\rightarrow
N_2
\rightarrow
N_3
\rightarrow
0
\]
be a short exact sequence of $\Gamma$-comodules.

For a morphism of $\Gamma$-comodules $h : X \rightarrow N_3$ and any lift
$\tilde{h}: X\rightarrow N_2$ of $h$ as a morphism of $A$-modules:
\begin{enumerate}
	\item
the morphism $\Delta \tilde{h}:= (\Gamma \otimes \tilde{h}) \circ \psi_X - \psi_{N_2} \circ \tilde{h}$
lies in the image of $\hom_A (X, \Gamma \otimes_A N_1 ) \hookrightarrow \hom_A (X, \Gamma \otimes_A N_2)$;
\item
the morphism $\Delta \tilde{h} : X \rightarrow \Gamma \otimes_A N_1$ is a cocycle and the associated class in $\ext^1 _\Gamma (X, N_1)$ satisfies
$[\Delta \tilde{h}] = \partial h$, where $\partial : \hom_\Gamma (X, N_3) \rightarrow \ext^1_\Gamma (X, N_1)$
is the connecting morphism.
\end{enumerate}
\end{lem}

In the presence of a given module splitting of a short exact sequence of comodules,
there is a useful description of a representing cocyle:

\begin{prop}
\label{prop:cocycle-canonical-splitting}
Let
\begin{eqnarray}
\label{eqn:ses-comod}
0 \rightarrow M
\rightarrow
\mathcal{E}
\stackrel{p}{\rightarrow}
X
\rightarrow
0
\end{eqnarray}
be a short exact sequence of $\Gamma$-comodules, in which $X$ is $A$-projective.

Then, for any choice of section $\sigma : X \rightarrow \mathcal{E}$, with associated
retract  $r :\mathcal{E}\rightarrow M$,  the composite:
\[
X \stackrel{- \sigma}{\rightarrow}
\mathcal{E}
\stackrel{\psi_{\mathcal{E}}}{\rightarrow}
\Gamma \otimes _A \mathcal{E}
\stackrel{\Gamma \otimes r}{\rightarrow}
\Gamma \otimes_A M
\]
is a cocycle which represents the extension class $\partial 1_X$, which corresponds to the short exact sequence (\ref{eqn:ses-comod}). A different choice of splitting as $A$-modules gives rise to a cohomologous cocycle.
\end{prop}

\begin{proof}
Consider the cocycle given by Lemma \ref{lem:explicit-connecting}, where $h$ is taken
to be the identity morphism of $X$ and $\tilde{h}$ the section $\sigma$. Then the
morphism $\Delta \tilde{h}$ is the difference $(\Gamma \otimes \sigma) \circ \psi_X - \psi_{\mathcal{E}} \circ \sigma$, considered (by abuse of notation)
also as a morphism to the image of $\Gamma \otimes_A M$ in $\Gamma \otimes_A \mathcal{E}$. Hence, composing with the retraction $\Gamma \otimes r$ gives $\Delta \tilde{h}=(\Gamma \otimes r) \Delta \tilde{h}$; now $r \circ \sigma =0$,  so that $\Delta \tilde{h} = -( \Gamma \otimes r) \circ \psi _{\mathcal{E}} \circ \sigma$, as required. The final statement is clear.
\end{proof}

\subsection{Base change and Landweber exactness}
\label{sect:cocycle-basechange}

Let $(A, \Gamma)$ be a flat Hopf algebroid and let $B
\stackrel{b}{\leftarrow} A \stackrel{c}{\rightarrow }C$ be ring morphisms.

\begin{nota}
	Write $\BGammaC$ for the ring $B \otimes _A \Gamma \otimes _A C$
 and  $\Gamma_B$ (respectively $\Gamma_C$) for the ring $ _B\Gamma_B$ (resp. $_C \Gamma _C$).
\end{nota}

There are induced Hopf algebroids $(B, \Gamma_B)$ and $(C, \Gamma_C)$.
	The ring $_B\Gamma_C$ has a left  $(B,\Gamma_B)$- right $(C,
\Gamma_C)$-bicomodule structure and the structure morphisms are
algebra morphisms.

\begin{defn}\cite{HS_comod}
\label{def:Landweber_exact}
A  morphism $b : A \rightarrow B$ is Landweber
exact with respect to $(A,\Gamma)$ if the functor $B\otimes_A - : \Gamma \dash \comod
\rightarrow B \dash \modules$ is exact.
\end{defn}

\begin{prop}
	\cite{HS_comod}
	If $b : A  \rightarrow B$ is Landweber exact with respect to $(A, \Gamma)$, then $(B, \Gamma_B)$ is a flat Hopf algebroid and $B \otimes_A -$ induces an exact functor $B\otimes_A - : \Gamma \dash \comod
\rightarrow \Gamma_B \dash \comod$.
\end{prop}

\begin{nota}
	For $f: X \rightarrow \Gamma \otimes_A M $ a morphism of left
$A$-modules, write $_Bf_C : B \otimes _A X \rightarrow \BGammaC \otimes _A M$ for the morphism of left $B$-modules:
\[
	B \otimes _A X
\stackrel{B \otimes f } {\rightarrow}
B \otimes _A \Gamma \otimes _A M
\stackrel{B \otimes A \otimes c \otimes M}{\rightarrow}
\BGammaC \otimes _A M.
\]
\end{nota}

\begin{exam}
 Let $X, M$ be as in Proposition \ref{prop:ext-cobar};	a cocycle $\kappa : X \rightarrow \Gamma \otimes_A M$  induces a
morphism $_B \kappa _C : B\otimes_A X \rightarrow \BGammaC \otimes _A M \cong \BGammaC \otimes_C (C \otimes_A M)$.
\end{exam}

\begin{lem}
\label{lem:cocycle_base_change}
Let $f : X \rightarrow \Gamma \otimes_A M $  and $g: M \rightarrow N$ be morphisms of left $A$-modules.
There is an identity
\[
	_B ((\Gamma \otimes_A g) \circ f ) _C = \Big( \BGammaC \otimes_C (C \otimes_A g) \Big) \circ \
_Bf_C.
\]
\end{lem}

\begin{lem}
	\label{lem:cocycle-base-change}
Let  $A\stackrel{b}{\rightarrow} B$ be
a Landweber exact morphism for $(A,\Gamma)$ and  $X,M$ be
$\Gamma$-comodules such that $X$ is  $A$-projective.

Let $\kappa : X \rightarrow \Gamma \otimes _A M$ be a cocycle representing an
extension $[\kappa] \in \ext^1_\Gamma (X, M)$. Then $B  \otimes_A X$ is a $B$-projective $\Gamma_B$-comodule and the morphism
$_B\kappa_B : B \otimes _A X \rightarrow \Gamma _B \otimes _B (B \otimes _A M)
$ is a cocycle which represents the class $B \otimes _A [\kappa ]   \in
\ext^1_{\Gamma_B} (B \otimes_A X , B \otimes_A M)$.
\end{lem}

\newcommand{\etalchar}[1]{$^{#1}$}
\providecommand{\bysame}{\leavevmode\hbox to3em{\hrulefill}\thinspace}
\providecommand{\MR}{\relax\ifhmode\unskip\space\fi MR }
\providecommand{\MRhref}[2]{%
  \href{http://www.ams.org/mathscinet-getitem?mr=#1}{#2}
}
\providecommand{\href}[2]{#2}

\end{document}